\documentclass[11pt]{amsart}
\usepackage{amsfonts,amssymb,amscd,amsmath,enumerate,verbatim,calc}

\catcode`@=11
\def\Proj{\operatorname{Proj}}  
\def\Spec{\operatorname{Spec}}
\def\height{\operatorname{ht}} 
\def\rdim{\operatorname{rdim}} 
\def\sdim{\operatorname{sdim}}

\def\sk{\par\smallskip}

\def\a{{\alpha}}
\def\b{{\beta}}
\def\l{{\lambda}}
\def\e{{\varepsilon}}
\def\Q{{\bf Q}}
\def\I{{\bf I}}
\def\J{{\bf J}}

\def\ps@icm{\def\@oddhead{\hfill \leftmark \hfill\thepage }
\def\@evenhead{\thepage \hfill \rightmark \hfill}
\def\@oddfoot+
\def\@evenfoot+}
\def\ps@first{\def\@oddhead{ICCM 2007 $\cdot$ Vol. II $\cdot$ 1--4
\hfill}
\def\@evenhead+
\def\@oddfoot+%\hfill \copyright\ China Higher Education Press}
\def\@evenfoot{\hfill\copyright\ China Higher Education Press}}
\def\list#1#2{\ifnum \@listdepth >5\relax \@toodeep \else \global
\advance \@listdepth\@ne \fi \rightmargin \z@ \listparindent\z@
\itemindent\z@ \csname @list\romannumeral\the\@listdepth\endcsname
\def\@itemlabel{#1}\let\makelabel\@mklab \@nmbrlistfalse #2\relax
\@trivlist \parskip -\parsep \parindent\listparindent \advance
\linewidth -\rightmargin \advance\linewidth -\leftmargin \advance
\@totalleftmargin \leftmargin \parshape \@ne \@totalleftmargin
\linewidth \ignorespaces}
\catcode`@=12
\def\thebibliography#1{*{References}
\list{[\arabic{enumi}]}{\settowidth \labelwidth{[#1]} \leftmargin
\labelwidth \advance \leftmargin \labelsep \usecounter{enumi}}
\def\newblock{\hskip .11em plus .33em minus .07em} \sloppy
\clubpenalty 4000 \widowpenalty 4000 \sfcode`\.=1000 \relax}

 \def\NN{{\mathbb N}} \def\ZZ{{\mathbb Z}}
\def\RR{{\mathbb R}} \def\CC{{\mathbb C}} 
 \def\PP{{\mathbb P}}

\def\frk{\mathfrak} \def\aa{{\frk a}} \def\pp{{\frk p}}
   
\def\mm{{\frk m}}  
\def\Phi{{\frk N}}
\def\Ac{{\mathcal A}}

\def\opn#1#2{\def#1{\operatorname{#2}}} % to make operators
\opn\chara{char} \opn\length{\ell} \opn\pd{pd} \opn\rk{rk}
\opn\projdim{proj\,dim} \opn\injdim{inj\,dim}
\opn\rank{rank} \opn\depth{depth} \opn\grade{grade} 
\opn\hei{ht} \opn\embdim{emb\,dim}\opn\codim{codim}
\opn\Tr{Tr} \opn\bigrank{big\,rank}
\opn\superheight{superheight} \opn\lcm{lcm}
\opn\rdim{rdim} \opn\trdeg{tr\,deg} \opn\reg{reg}  \opn\lreg{lreg} 
\opn\ini{in} \opn\lpd{lpd} \opn\size{size} \opn{\mult}{mult}
\opn\div{div} \opn\Div{Div} \opn\cl{cl} \opn\Cl{Cl}
\opn\Spec{Spec} \opn\Supp{Supp} \opn\supp{supp} 
\opn\Sing{Sing} \opn\Ass{Ass} \opn\Min{Min}
\opn\Proj{Proj}
\opn\Ann{Ann} \opn\Rad{Rad} \opn\Soc{Soc}
\opn\Syz{Syz} \opn\Im{Im} \opn\Ker{Ker} \opn\Coker{Coker}
\opn\Am{Am} \opn\Hom{Hom} \opn\Tor{Tor} \opn\Ext{Ext}
\opn\End{End} \opn\Aut{Aut} \opn\id{id}

\opn\nat{nat} \opn\pff{pf} % \pf exists already
\opn\Pf{Pf} \opn\GL{GL} \opn\SL{SL} \opn\mod{mod} \opn\ord{ord}
\opn\Gin{Gin} \opn\Hilb{Hilb}
\opn\adeg{adeg} \opn\std{std}\opn\ip{infpt}
\opn\Pol{Pol} \opn\sat{sat} \opn\Var{Var}
\opn\aff{aff} \opn\con{conv} \opn\relint{relint} \opn\st{st}
\opn\lk{lk} \opn\cn{cn} \opn\core{core} \opn\vol{vol}
\opn\link{link} \opn\star{star}
\opn\gr{gr}

\def\pot#1#2{#1[\kern-0.28ex[#2]\kern-0.28ex]}

\opn\dirlim{\underrightarrow{\lim}}
\opn\inivlim{\underleftarrow{\lim}}

\newtheorem{Theorem}{Theorem}[section]
\newtheorem{Lemma}[Theorem]{Lemma}
\newtheorem{Corollary}[Theorem]{Corollary}
\newtheorem{Proposition}[Theorem]{Proposition}

\newtheorem{Example}[Theorem]{Example}

\newtheorem{Question}[Theorem]{Question}
\let\epsilon\varepsilon
\let\phi=\varphi
\let\kappa=\varkappa
\let\lm=\lambda

\opn\dis{dis}
\opn\Lex{Lex}
\begin{document}
\title[Hilbert functions and mixed multiplicities]
{Hilbert functions of multigraded algebras,  mixed 
multiplicities of ideals and their applications}
\author{N. V. Trung}
\address{N. V. Trung, Institute of Mathematics, 18 Hoang Quoc Viet, 10307, Hanoi, Vietnam}
\email{nvtrung@math.ac.vn}
\author{J. K. Verma}
\address{J. K. Verma, Department of Mathematics, Indian Institute
of Technology Bombay, Mumbai, India}
\email{jkv@math.iitb.ac.in}

\begin{abstract} 
This paper is a survey on major results on Hilbert functions of multigraded algebras and mixed multiplicities of ideals, 
including their applications to the computation of Milnor numbers of complex analytic hypersurfaces with isolated singularity, multiplicities of blowup algebras and  mixed volumes of polytopes.  
\end{abstract}

\maketitle
\thispagestyle{empty}

\tableofcontents
\section{Introduction}

Let $R=\bigoplus_{t=0}^\infty R_t$ be a 
Noetherian  graded algebra over a field $k$. Then $R_t$ is a finite dimensional 
$k$-vector space. In his historic paper, \cite{hi1} Hilbert considered 
the generating function 
$$H(R,z):=\sum_{t=0}^\infty H_R(t)z^t$$ 
of the sequence $H_R(t):=\dim_kR_t.$ By using his {\em Syzygy Theorem}, he proved that
if $R=k[f_1,f_2,\ldots,f_s]$ where $f_i \in R_{d_i}$ for $i=1,2,\ldots,s,$
then there exists a polynomial $h(z) \in \ZZ[z]$ such that
$$ H(R,z)=\frac{h(z)}{(1-z^{d_1})(1-z^{d_2})\cdots(1-z^{d_s})}.$$

We say that $R$ is standard if $R$ is generated over $R_0$ by elements of degree 1. 
In this case
the Hilbert function $H_R(t)$ is given by a polynomial $P_R(t)$ for all $t$ 
large enough.  Lasker \cite{la} showed that the Krull dimension of 
$R,$ denoted by $\dim R,$ is $\deg P_R(t) +1.$ In the same paper, Lasker indicated that  these results
could be generalized to Hilbert functions of $\NN^r$-graded algebras.
The ideas  of Lasker and Noether were  presented by Van der Waerden in 
a detailed exposition \cite{w}. 

Let $X$ and $Y$ be two sets of $m+1$ and $n+1$ indeterminates, respectively.
A polynomial in $k[X,Y]$ is bihomogeneous if it is  homogeneous in $X$ and $Y$ separately.
An ideal $I$ is called bihomogeneous if it is generated by bihomogeneous polynomials. 
Let $V \subset \PP^m \times \PP^n$ be the zero set of a collection of bihomogeneous polynomials.
Then the ideal $I(V)$  of polynomials in $k[X,Y]$ which vanish on $V$ is 
bihomogeneous. Therefore, the coordinate ring $k[X,Y]/I(V)$ is a bigraded algebra of the form 
$$R = \bigoplus_{(u,v)\in \NN^2}R_{(u,v)},$$
where $R_{(u,v)}$ is a finite dimensional $k$-vector space.

Van der Waerden showed that $H_R(u,v) := \dim_k R_{(u,v)}$ is given by a polynomial $P_R(u,v)$
with rational coefficients for all large values of $u,v.$ The degree of $P_R(u,v)$
is at most $\dim R-2.$ Let $r=\deg P_R(u,v).$ Write $P_R(u,v)$ in the form 
$$
P_R(u,v)=\sum_{i+j \leq r} e_{ij}(R) \binom{u}{i}\binom{v}{j}.
$$
We call $P_R(u,v)$ the {\it Hilbert polynomial} and the numbers $e_{ij}(R)$ with $i+j = r$ the {\it mixed multiplicities} of $R$. The  mixed multiplicities have geometrical significance.\medskip

\noindent {\bf Theorem} (Van der Waerden, 1928) {\em Let $P$ be a bihomogeneous prime ideal
of $k[X,Y]$ and $R = k[X,Y]/P$. Then  $e_{ij}(R)$ is the number of points of intersection of the variety 
$$V(P) = \{\a \in \PP^m\times \PP^n|\ f(\a) = 0\ \forall\  f \in P\}$$ 
with a linear space defined by $i$ general linear equations in $X$ and $j$ general linear equations 
in $Y.$ }

\medskip

The Hilbert polynomial $P_R(u,v)$ and the mixed multiplicities $e_{ij}(R)$ can be defined for any Noetherian bigraded algebra $R$ over an Artinian local ring which is standard in the sense that it is generated by elements of degree $(1,0)$ and $(0,1)$. These objects  were not so well studied until recently. The total degree of $P_R(u,v)$ was characterized independently by Schiffel  \cite{sc} and Verma-Katz-Mandal \cite{vkm}. They showed that $\deg P_R(u,v)+2$ is the maximal dimension of the relevant prime ideals of $R.$ 
Verma-Katz-Mandal also showed that the mixed multiplicities $e_{ij}(R)$ can be any sequence of non-negative integers with at least a positive entry. Recently, Trung was able to characterize the  degrees in $u$ and $v$ of $P_R(u,v)$ and the positive mixed multiplicities in \cite{Tr1}. In particular, he showed that the range of the positive mixed multiplicities is rigid if $R$ is a domain or a Cohen-Macaulay ring, thereby solving an open question of Verma-Katz-Mandal. 

An important case of mixed multiplicities of a bigraded algebra is  the mixed multiplicities of two ideals.
Let $(A,\mm)$ be a local ring. For any pair of  $\mm$-primary ideals $I$ and $J,$ one can consider the length function $\ell(A/I^uJ^v)$ which is the sum transform of the Hilbert function of the standard bigraded algebra
$$R(I|J) := \bigoplus_{u,v \ge 0} I^uJ^v/I^{u+1}J^v$$ 
over the quotient ring $A/I$. Bhattacharya \cite{bh} showed 
that this function is given by a polynomial $P(u,v)$ of degree $d =\dim A$ and that  it can be written as 
$$P(u,v)=\sum_{i+j\le d}a_{ij}(I|J)\binom{u+i}{i}\binom{v+j}{j}$$
for certain  integers $a_{ij}(I|J).$ 
We set $e_j(I|J):= a_{ij}(I|J)$ for $i+j= d$.   
These integers were named later as mixed multiplicities by Teissier in \cite{t1} where 
he found significant applications of $e_j(I|J)$ in the study of singularities of complex analytic hypersurfaces.
In particular, the Milnor numbers  
of linear sections of a complex analytic hypersurface at an isolated singularity are exactly the mixed multiplicities of the maximal ideal and the Jacobian ideal of the hypersurface.

Teissier found several interesting properties of mixed multiplicities of ideals
which have inspired subsequent works substantially.
His characterization of mixed multiplicities as Samuel's multiplicities of general elements led Rees \cite{r4} to the introduction of joint reductions of ideals which generalize the important concept of reduction of an ideal  in  multiplicity theory. 
Another instance is the inequalities 
$$e_j(I|J)^d \le e(I)^{d-j}e(J)^j$$
for $j = 0,...,d$ which implies the Minkowski inequality
$$e(IJ)^{1/d} \le e(I)^{1/d}+e(J)^{1/d}.$$
Teissier raised it as a conjecture and showed   it for reduced Cohen-Macaulay complex analytic algebras \cite{t3}.
This conjecture was settled in the affirmative by Rees and Sharp in \cite{resh}.

Mixed multiplicities are also  defined  if $I$ is an $\mm$-primary ideal and $J$ an arbitrary ideal by using the standard graded algebra $R(I|J)$. Katz-Verma \cite{kv} and Verma  \cite{v2} \cite{v7} studied first the mixed multiplicities in these cases. They showed that these mixed multiplicities can be used to compute the multiplicity of the Rees algebra and the extended Rees algebra. D'Cruz \cite{d3} obtained multiplicity formula for multigraded extended Rees algebra. Herzog-Trung-Ulrich \cite{htu} have  devised an effective  method to compute the multiplicity of the Rees algebras which is similar to that of Gr\"obner bases. This method has been exploited by 
Hoang \cite{Ho} \cite{Ho2} and   Raghavan-Verma \cite{rv} to compute mixed multiplicities of ideals generated by $d$-sequences, quadratic sequences and filter-regular sequences of homogeneous elements of non-decreasing degrees.

A systematic study of mixed multiplicities of two ideals in the general case was carried out by Trung in \cite{Tr1}. He characterized the positive mixed multiplicities  and showed how to compute them by means of  general elements. As a consequence, the range of the positive mixed multiplicities is rigid and depends only on the ideal $J$.

Mixed multiplicities are  also defined for an $\mm$-primary ideal and a sequence of ideals of $A$.
To handle the complexity of this case Trung and Verma \cite{TV} used a multigraded version of the associated graded ring in order to introduce 
superficial sequence for a collection of ideals. Using this notion they obtained  similar results as in the case of two ideals.
These results can be applied to describe mixed volumes of lattice polytopes as mixed multiplicities, thereby giving a purely algebraic proof of Bernstein's theorem which asserts that the number of common zeros of a  system of $n$ Laurent polynomials in $n$ indeterminates with 
finitely many   zeroes in the torus  $(\CC^*)^n$ is bounded above by the mixed volume of their Newton polytopes. 
  
Another interesting instance of mixed multiplicities is the multiplicity sequence of an ideal introduced by Achilles and Manaresi in \cite{AM2}.
Let $I$ be an arbitrary ideal in  a local ring $(A,\mm)$. The associated graded ring
$$R  = \bigoplus_{(u,v) \in \NN^2}\big(\mm^uI^v + I^{v+1}/\mm^{u+1}I^v +I^{v+1}\big),$$
 is a standard bigraded algebra over the residue field $A/\mm$. 
The sum transform
$$H_R^{(1,1)}(u,v) := \sum_{i = 0}^u\sum_{j = 0}^vH_R(i,j)$$
of the Hilbert function of $R$ is given by a polynomial $P_R^{(1,1)}(u,v)$ of degree $d$ for $u,v$ large enough. 
If we write this polynomial in the form
$$P_R^{(1,1)}(u,v) =\sum_{i= 0}^d \frac{c_{ij}(R)}{i!(d-i)!}u^iv^{d-i}  +
\text{\rm lower-degree terms},$$
then $c_i(I) := c_{i\; d-i}(R) $ are non-negative integers for $i = 0,...,d$. 
Achilles and Manaresi  call $c_0(I),...,c_d(I)$ the {\it multiplicity sequence} of $I$.
The multiplicity sequence can be considered as a generalization of the multiplicity of an $\mm$-primary ideal. 
In fact, if $I$ is an $\mm$-primary ideal, then $c_0(I) = e(I)$ and $c_i(I) = 0$ for $i > 0$.
In particular, $c_0(I) > 0$ if and only if the analytic spread of $I$, $s(I),$ equals $d.$  In this case,  $ c_0(I)$ is called the $j$-{\it multiplicity} of $I$ \cite{AM1}. Flenner-Manaresi \cite{FM} used $j$-multiplicity to give a numerical criterion for reduction of ideals.
The multiplicity sequence can be computed by means of the intersection algorithm which was introduced by St\"uckrad-Vogel \cite{SV} in order to prove a refined version of the Bezout's theorem.

In general, the Hilbert function $H_R(u,v)$ of a finitely generated bigraded algebra $R$ over a field $k$ is not a polynomial for large $u,v$. However, if $R$ is generated by elements of 
bidegrees $(1,0),(d_1,1),\ldots,(d_r,1)$, where $d_1,\ldots,d_r$ are non-negative integers, then 
there exist integers $c$ and $v_0$ such that $H_R(u,v)$ 
is equal to a polynomial $P_R(u,v)$ for $u \ge cv$ and $v \ge v_0$. 
This case was considered first by P.~Roberts in \cite{ro} and then by Hoang-Trung in \cite{HT}.  
Hilbert polynomials of bigraded algebras of 
the above type appear in Gabber's proof of Serre's non-negativity 
conjecture \cite{Ro2} for intersection multiplicities and that the positivity of certain coefficient 
of such a Hilbert polynomial is strongly related to Serre's positivity
conjecture on intersection multiplicities \cite{Ro3}. 

An instance of such non-standard bigraded algebra is the Rees algebra of a homogeneous ideal $I$ in a standard graded algebra $A$. The existence of the Hilbert polynomial in this case allows us to study the behavior of the Hilbert polynomials of the quotient rings $A/I^v$ for $v$ large enough \cite{HPV}. We can also use the mixed multiplicities of the Rees algebra of $I$ to compute the degree of the embedded varieties of the blow-ups of $\Proj A$ along $I$.

The above development of the theory of Hilbert functions of multigraded algebras and of mixed multiplicities of ideals will be discussed in more detail in subsequent sections.

%1. Hilbert functions of multigraded algebras
%\par 2. Positivity of mixed multiplicities
%\par 3. Mixed multiplicities of ideals: the $\mm$-primary case
%\par 4. Mixed multiplicities of two ideals: the general case
%\par 5. numbers and mixed multiplicities
%\par 6. Multiplicities of blow-up algebras
%\par 7. Mixed  multiplicities of a sequence of ideals: the general case
%\par 8. Mixed volume of lattice polytopes
%\par 9. Minkowski inequalities and equalities
%\par 10. The multiplicity sequence
%\par 11. Hilbert function of non-standard bigraded algebras
%\par 12. Hilbert function of bigraded Rees %algebras
%\smallskip

Illustrating examples and open problems for  further study will be given throughout the paper.
The  results discussed in this paper merely reflect  authors' interests and do not cover all the developments due to  lack of space and time and also due to  their ignorance.

\noindent
{\bf Acknowledgements:} The authors thank Bernard Teissier
for several useful comments.

\section{Hilbert functions of multigraded algebras}

Let $R = \bigoplus_{(u,v) \in \NN^2} R_{(u,v)}$ be a Noetherian bigraded algebra over an Artinian local ring $R_{(0,0)}=k.$
We define the Hilbert function of $R$ by  $H_R(u,v) := \ell(R_{(u,v)})$,
where $\ell$ denotes the length.

If  $R$ is standard graded,  $H_R(u,v)$ is given by a polynomial $P_R(u,v)$ for all $u, v$ large enough.
In order to determine the total degree of $P_R(u,v)$ we need the following notions. 

We say that a bihomogeneous ideal $I$ of $R$
is called {\em irrelevant} if $I_{(u,v)}=R_{(u,v)}$ for all  $u,v$ large. We say that $I$ is 
{\em relevant}  if it is not irrelevant. Let $\Proj R$ denote the set of all bihomogeneous relevant prime ideals of $R$.  
The {\it relevant dimension} $\rdim R$ of 
$R$ is defined by 
$$ \rdim R=\max \{\dim R/P|\ P \in \Proj R \}.$$ 
The total degree of $P_R(u,v)$ was found independently by Schiffel \cite{sc} and Katz-Mandal-Verma \cite{vkm}.

\begin{Theorem}
$\deg P_R(u,v) =\rdim R-2.$
\end{Theorem}

Let $r = \rdim R-2$. As in the case $k$ is a field, if we write $P_R(u,v)$ in the form 
$$
P_R(u,v)=\sum_{i+j \leq r} e_{ij}(R) \binom{u+i}{i}\binom{v+j}{j},
$$
then the numbers $e_{ij}(R)$  are non-negative integers for all $i,j$ with $i+j = r$.
These numbers are called the {\it mixed multiplicities} of $R$. 

\begin{Example}
{\rm Let $R = k[x_1,...,x_m,y_1,...,y_n]$ with $\deg x_i = (1,0)$ and $\deg y_j = (0,1)$.
Then
$$H_R(u,v)= P_R(u,v) = \binom{u+m-1}{m-1}\binom{v+n-1}{n-1}$$
for all $(u,v) \in \NN^2$. Therefore, $\deg P_R(u,v) = m+n-2$ and} 
$$e_{ij}(R) = \left\{\begin{array}{ll} 1 &\ \text{ if }\ i = m-1, j = n-1,\\
0 & \text{\;\;\;otherwise}.\end{array} \right.$$
\end{Example}

The computation of mixed multiplicities can be reduced to the case of a bigraded domain by using the following associativity formula
(see e.g. \cite{HT}).

\begin{Proposition} \label{associative}
Let ${\mathcal A}(R)$ be the set of the prime ideals $P \in 
\Proj R$ with $\dim R/P = \rdim R$. Then
$$e_{ij}(R)= \sum_{P \in {\mathcal A}(R)}\ell(R_P)e_{ij}(R/P).$$
\end{Proposition}

Katz, Mandal and Verma \cite{vkm} showed that the mixed multiplicities can be any sequence of non-negative integers with at least a positive entry. 

\begin{Example}
{\rm Let $a_0,...,a_n$ be an arbitrary sequence of non-negative integers with at least a positive entry.
Let $S = k[x_0,...,x_n,y_0,y_1,...,y_n]$ be a bigraded polynomial ring with $\deg x_i = (1,0)$ and $\deg y_j = (0,1)$. Let $Q_t  = (x_0^{a_t},x_1,...,x_{n-t-1},y_0,y_1,...,y_{t-1})$, $t = 0,...,n$.
Set $R = S/Q_0 \cap Q_1 \cap \cdots \cap Q_n$. Then $\rdim R = n$ and
${\mathcal A}(R) = \{P_0,...,P_n\}$, where
$P_t = (x_0,x_1,...,x_{n-t-1},y_1,...,y_{t-1})$. We have $\ell(R_{P_t}) = a_t$ and
$$e_{in-i}(R/P_t) = \left\{\begin{array}{ll} 1 &\ \text{ if }\ i = t,\\
0 & \text{otherwise}.\end{array} \right.$$
By the associativity formula we get 
$e_{in-i} = a_i$ for $i = 0,...,n$.}
\end{Example}

A standard way to make a standard  bigraded  algebra $R$ into an 
$\NN$-graded algebra is by defining $R_t =\bigoplus_{u+v=t} R_{(u,v)}.$ 
This algebra is obviously standard graded.

For any Noetherian standard graded $\NN$-graded algebra $R$ over an Artinian local ring, we have
$\rdim R = \dim R$. Let $d = \dim R$. If we write 
$P_R(t)$ in the form
$$P_R(t)=\sum_{i=0}^{d-1} a_i(R) \binom{t+d-1-i}{d-1-i},$$
then $e(R) := a_0(R)$ is called the {\it multiplicity} of $R$.

The relationship between multiplicity and mixed multiplicities  was found independently in the unpublished
thesis of Dade \cite{da} and in \cite{vkm}.

\begin{Theorem}[Dade, 1960; Katz-Mandal-Verma, 1994] \label{total}
Let $R$ be a Noetherian bigraded algebra over an Artinian local ring.
Assume that the ideals $(R_{(1,0)})$ and $(R_{(0,1)})$  have positive height. Then
$$e(R)=\sum_{i+j=\dim R-2} e_{ij}(R).$$
\end{Theorem}

The above results for bigraded algebras have been extended to multigraded modules by Herrmann-Hyry-Ribbe-Tang in the paper \cite{hhrt}. To summarize their results we fix some notations. 

Let $s$ be any non-negative integer. Let $R = \oplus_{u \in \NN^s}R_u$ 
be a  Noetherian standard $\NN^s$-graded algebra 
over an Artinian local ring $k$, where \lq standard\rq\ means
$R$ is  generated by homogeneous  elements of degrees $(0,..,1,..,0)$, where $1$ 
occurs only as the $i$th component, 
$i = 1,...,s.$

For $\a = (\a_1,...,\a_s)$ and $\b = (\b_1,...,\b_s)$ we write $\a > \b$ 
if $\a_i > \b_i$  for all $i=1,\ldots, s.$ 
Let 
$$R_+ = \bigoplus_{ \a > {\bf 0}} R_\a.$$
We define $\Proj R$  to be the set of all $\NN^s$-graded prime ideals which do not contain 
$R_+$. It is easy to see that $P \in \Proj R$ if and only if  $P_u \neq R_u$ for all $u \in \NN^s$.

Let $M = \bigoplus_{u \in \ZZ^s}M_u$ be a finitely generated $\ZZ^s$-graded module over $R$. 
Then $M_u$ is a $k$-module of finite length. We call $H_M(u) := \ell(M_u)$  the
{\it Hilbert function} of $M$. Moreover, we define the {\it relevant dimension} of $R$ to be the number
$$\rdim M := \max\{\dim R/P|\ P \in \Proj R\ \text{and}\ M_P \neq 0\}.$$

\begin{Theorem}[Herrmann-Hyry-Ribbe-Tang, 1997] 
For $u \gg 0$, $H_M(u)$ is given by a polynomial $P_M(u)$ with rational coefficients
having total degree  $\rdim M-s$.
\end{Theorem}

Let $r = \rdim R-s$. If we write $P_M(u)$ in the form
$$P_M(u) = \sum_{\a \in \NN^s, |\a| = r} \frac{1}{\a!}e_\a(M)u^\a  +
\text{terms of degree $< r$},$$
where $\a = (\a_1,...,\a_s)$ with 
$$
|\a| := \a_1 + \cdots + \a_s, \ \
\a!  := \a_1!\cdots\a_s!,\ \ \mbox{and} \ \
u^\a  := u_1^{\a_1}\cdots u_s^{\a_s},
$$
then $e_\a(M)$ are non-negative integers if $|\a| = r$. We call these coefficients  the {\it mixed multiplicities} of 
$M$.  \par

Now we  consider the difference of the Hilbert function and the Hilbert polynomial.

Let  $R$ be a Noetherian standard $\NN$-graded algebra over an Artinian local ring.
Let $M=\bigoplus_{t\in \ZZ} M_t$ be  a finite $\ZZ$-graded  module over $R$. 
Let $H^i_{R_+}(M)$ denote the $i$th local cohomology module of $M$ with respect to $R_+$.
Then for all $t \in \ZZ$,
$$H_M(t)-P_M(t)=\sum_{i\ge 0}(-1)^i\ell(H^i_{R_+}(M)_t),$$ 
which is known as the Grothendieck-Serre formula. 

A similar formula for the bigraded case was proved by Jayanthan and Verma \cite{jv}.

\begin{Theorem}[Jayanthan-Verma, 2002]
Let $R$ be a Noetherian standard bigraded algebra over an Artinian local ring. 
Let $M$ be a finite bigraded module over $R$.  Then for all $u,v \in \ZZ$,
$$
H_M(u,v) - P_M(u,v) = \sum_{i\geq 0}(-1)^i\ell(H^i_{R_+}(M)_{(u,v)}).
$$
\end{Theorem}

\section{Positivity of mixed multiplicities}

Let $R$ be a Noetherian standard bigraded algebra over an Artinian local ring.
Can we say which mixed multiplicities of $R$ 
 are  positive ?
To give an answer to this question let us first express the total degree of the Hilbert polynomial $P_R(u,v)$ in another way. For any pair of ideals $\mathfrak a, \mathfrak b$ let
$${\mathfrak a}:{\mathfrak b}^\infty := \{x \in S|\ \text{there is a positive integer $n$ such that $x{\mathfrak b}^n \subseteq {\mathfrak a}$}\}.$$ 
It is easy to see that  $\rdim R = \dim R/0:R_+^\infty$ and therefore
$$\deg P_R(u,v) = \dim R/0:R_+^\infty-2.$$

Similarly, one can also compute the partial degrees of the Hilbert polynomial $P_R(u,v)$ \cite{Tr2}.

\begin{Theorem}[Trung, 2001] \label{degree} 
\begin{eqnarray*}  \deg_uP_R(u,v) & = & \dim R/(0:R_+^\infty+(R_{(0,1)}))-1,\\ 
\deg_vP_R(u,v) & = & \dim R/(0:R_+^\infty+(R_{(1,0)}))-1. \end{eqnarray*}
\end{Theorem}

For simplicity we set 
\begin{eqnarray*} 
r & := & \dim R/0:R_+^\infty-2,\\ 
r_1 & := & \dim R/(0:R_+^\infty+(R_{(0,1)}))-1,\\ 
r_2 & := & \dim R/(0:R_+^\infty + (R_{(1,0)}))-1. 
\end{eqnarray*}

\begin{Corollary}\label{zero} 
$e_{ij}(R) = 0$ for $i > r_1$ or $j > r_2$. 
\end{Corollary}

To characterize the positive mixed multiplicities we shall need the concept of a filter regular sequence which  originates from the theory of generalized Cohen-Macaulay rings \cite{CST}. A sequence $z_1,\ldots,z_s$ of homogeneous elements in $R$ is called {\it filter-regular}  if for $i = 1,\ldots,s$, we have 
$$[(z_1,\ldots,z_{i-1}):z_i]_{(u,v)} = (z_1,\ldots,z_{i-1})_{(u,v)}$$
for $u$ and $v$ large enough. It is easy to see that $z_1,\ldots,z_s$ is filter-regular  if and only if $z_i \not\in P$ for all associated prime ideals $P \not\supseteq R_+$ of $(z_1,\ldots,z_{i-1})$, $i = 1,\ldots,s$  (see e.g. \cite{Tr} for basic properties). 

The following result gives an effective criterion for the positivity of a mixed multiplicity $e_{ij}(R)$ and shows how to compute $e_{ij}(R)$ as the multiplicity of an $\NN$-graded algebra \cite{Tr2}.

\begin{Theorem}[Trung, 2001]  \label{nonzero} 
Let $i, j$ be non-negative integers with $i+j = r$. 
Let $x_1,\ldots,x_i$ be a filter-regular sequence of homogeneous elements of degree $(1,0)$. Then $e_{ij}(R) > 0$ if and only if 
$$\dim R/((x_1,\ldots,x_i):R_+^\infty+(R_{(0,1)})) = j+1.$$ 
In this case, if we choose homogeneous elements $y_1,\ldots,y_j$ of degree $(0,1)$ such that $x_1,\ldots,x_i,y_1,\ldots,y_j$ is a filter-regular sequence, then  
$$e_{ij}(R) = e(R/(x_1,\ldots,x_i,y_1,\ldots,y_j):R_+^\infty).$$ 
\end{Theorem}

If the residue field of $R_0$ is infinite, one can always find homogeneous elements $x_1,\ldots,x_i$ of degrees $(1,0)$ and $y_1,\ldots,y_j$  of degree $(0,1)$ such that $x_1,\ldots,x_i,y_1,\ldots,y_j$ is a filter-regular sequence. 

For $i = 0$ we get the condition $r_2 +1 = r$ which yields the following criterion for the positivity of $e_{0r}$.

\begin{Corollary} \label{first}  
$e_{0r}(R) > 0$ $(e_{r0}(R) > 0)$ if and only if $r_2+1 = r$ $(r_1+1 = r)$. 
\end{Corollary} 

In spite of Corollary \ref{zero} one might ask whether $e_{i\;r-i}(R) > 0$ for $i = r_1,r-r_2$. 
Using Theorem \ref{nonzero} one can easily construct examples with $e_{i\;r-i}(R) = 0$ for $i = r_1,r-r_2$.

\begin{Example} 
{\rm Let $R = k[X,Y]/(x_1,y_1) \cap (x_1,x_2,x_3) \cap (y_1,y_2,y_3)$ with $X = \{x_1,x_2,x_3,x_4\}$, $Y = \{y_1,y_2,y_3,y_4\}$ and $\deg x_i = (1,0),\ \deg y_i = (0,1),\ i = 1,2,3,4$. Then $R/(R_{(1,0)}) \cong k[Y]$ and $R/(R_{(0,1)}) \cong k[X]$. Since $0:R_+^\infty = 0$, we get 
\begin{align*}
r & = \dim R - 2 = 4,\\
r_1 &= \dim R/(R_{(1,0)}) - 1= 3,\\
r_2 & = \dim R/(R_{(0,1)}) - 1 = 3.
\end{align*}
It is clear that $x_4$ is a non-zerodivisor in $R$. Since 
$x_4R:R_+^\infty + (R_{(1,0)}) = (x_1,x_2,x_3,x_4,y_1)R,$  
we have 
$$\dim R/(x_4R:R_+^\infty + (R_{(1,0)}) = \dim k[y_2,y_3,y_4] = 3 < 3+1.$$
Hence $e_{13}(R) = 0$. By symmetry we also have $e_{31}(R) = 0$. Now we want to compute the only non-vanishing mixed multiplicity $e_{22}(R)$ of $R$. It is easy to check that $x_4,x_2,y_4,y_2$ is a filter-regular sequence in $R$.
Put $Q = (x_4,x_2,y_4,y_2)$. Then
$$R/Q:R_+^\infty = k[X,Y]/(x_1,x_2,x_4,y_1,y_2,y_4) \cong k[x_3,y_3].$$ 
Hence $e_{22}(R) = e(R/Q:R_+^\infty) = \ell(k).$}
\end{Example} 

We say that the sequence of positive mixed multiplicities is {\it rigid} if there are integers $a, b$ such that $e_{i\;r-i}(R) > 0$ for $a \le i \le b$ and $e_{i\;r-i}(R) = 0$ otherwise. Obviously, that is the case if $e_{i\;r-i}(R) > 0$ for $r-r_2 \le i \le r_1$.

Katz, Mandal and Verma \cite{vkm} raised the question whether the sequence of positive mixed multiplicities is rigid if $R$ is a domain or Cohen-Macaulay. We shall see that this question has a positive answer by showing that $e_{i\;r-i}(R) > 0$ for $r-r_2 \le i \le r_1$ in these cases.

Recall that a commutative noetherian ring $S$ is said to be {\it connected in codimension 1} if the minimal dimension of closed subsets 
$Z \subseteq \Spec(S)$ for which  $\Spec(S)\setminus Z$ is disconnected is equal to $\dim S-1$.

Using a version of Grothendieck's Connectedness Theorem  due to Brodmann and Rung \cite{br} one can prove the following sufficient condition for  the rigidity of mixed multiplicities.

\begin{Theorem}[Trung, 2001] \label{rigid1} 
Assume that all maximal chains of prime ideals in $R/0:R_+^\infty$ have the same length.  
Then $e_{i\;r-i}(R)> 0$ for $i = r-r_2, r_1$. If $R/0:R_+^\infty$ is moreover connected in codimension 1, then  $e_{i\;r-i}(R)> 0$ for $r-r_2 \le i \le r_1$. 
\end{Theorem}

If $R$ is a domain or a Cohen-Macaulay ring with $\height R_+ \ge 1$, then $R$ is connected in codimension 1 by
Hartshorne's Connectedness Theorem. Hence the sequence of positive mixed multiplicities is rigid in these cases.

\begin{Corollary} \label{rigid2} 
Let $R$ be a domain or a Cohen-Macaulay ring with $\height R_+ \ge 1$. Then $e_{i\;r-i}(R)> 0$ for $r-r_2 \le i \le r_1$. 
\end{Corollary}

\section{Mixed multiplicities of ideals: the $\mm$-primary case} 

Throughout this section $(A, \mm)$ will denote  Noetherian local ring of positive dimension $d$ with infinite residue field. 

Let $I$ be an $\mm$-primary ideals. Then $A/I^t$ is of finite length for all $t \ge 0$. 
It is well-known that the  function $\ell (A/I^t)$ is given by a polynomial $P_I(t)$ for large $t,$ called the {\it Hilbert-Samuel polynomial} of $I.$ The degree of $P_I(t)$ is $d$ and if we write   $P_I(t)$ in terms of binomial coefficients as:
$$P_I(t)=e_0(I)\binom{t+d-1}{d}-e_1(I)\binom{t+d-2}{d-1}+\cdots+(-1)^de_d(I).$$

The coefficients $e_0(I), e_1(I),\ldots, e_d(I)$ are integers. 
The coefficient  $e_0(I)$ is a positive integer called the  ({\it Samuel's}) {\it multiplicity} of $I$  and it will be denoted by  
$e(I).$

Let $J$ be an $\mm$-primary ideal (not necessarily different from $I$). Bhattacharya \cite{bh} showed a similar property for  the bivariate function 
$$B(u, v) := \ell(A/I^uJ^v).$$

\begin{Theorem}[Bhattacharya, 1955] 
There exists a polynomial $P(u,v)$ of total degree $d$ in $u$ and $v$ with rational coefficients so that $B(u,v) = P (u,v)$ for all large $u,v.$ 
The terms of total degree $d$ in $P(u,v)$ have the form
$$\frac{1} {d!}\left\{e_0(I|J)r^d  + \cdots  +\binom{d}{i}e_i(I|J)u^{d-i}v^i + \cdots + e_d(I|J)s^d\right\}$$
where  $e_0(I|J),...,e_d(I|J)$ are certain positive
integers.
\end{Theorem}
The numbers $e_0(I|J ), \ldots , e_i(I|J), \ldots, e_d(I|J)$ 
were termed as the {\it mixed multiplicities}  of $I$ and $J$ by Teissier \cite{t1}. 
We have  the following relationship between mixed multiplicities and multiplicity.

\begin{Proposition} [Rees, 1961] 
$e_0(I|J) = e(I)$ and  $e_d(I|J) = e(J).$ 
\end{Proposition}

For all positive numbers $u,v$ we have  $P_{I^uJ^v}(t) = P(ut,vt)$. Therefore,
$$
e(I^uJ^v) = e_0(I|J)u^d  +\cdots  +\binom{d}{i}e_i(I|J)u^{d-i}v^i  +\cdots + e_d(I|J)v^d.
$$
This is perhaps the reason, why Teissier defined
$e_0(I|J), \ldots, e_d(I|J),$ to be the mixed multiplicities of the ideals $I$ and $J.$

We shall see later that mixed multiplicities can  always be expressed as Samuel's multiplicities.
There are numerous ways of computing multiplicity of an $\mm$-primary ideal.
We refer the reader to \cite[Chapter 11]{sh}. In particular, if $A$ is a Cohen-Macaulay ring and $I$ is a parameter ideal, then $e(I) = \ell(A/I)$.

An effective way for the computation of multiplicity was discovered by Northcott and Rees \cite{nr} by using the following notion: An ideal $J \subseteq I$ is called a {\it reduction} of $I$ if there exists an integer $n$ such that $JI^n = I^{n+1}.$ They showed that any minimal reduction of an $\mm$-primary ideal $I$ is a parameter ideal if $R/\mm$ is infinite.  It is easy to check that if  $J$ is a reduction of $I$ then  $e(I)=e(J).$

There is a deep connection between reductions and multiplicity \cite{r1}. 

\begin{Theorem} [Rees Multiplicity Theorem, 1961] 
Let $A$ be a quasi-unmixed local ring. Let $J \subseteq I$ be  $\mm$-primary ideals.  Then $J$ is a reduction of $I$ if and only if  $e(I)=e(J).$
\end{Theorem}

Recall that $A$ is called a {\em quasi-unmixed local ring} if 
$\dim \hat{A}/p=\dim R$ for each minimal prime $p$ of $\hat{A}$ where
the $\mm$-adic completion of $A$ is denoted by $\hat{A}.$

To generalize Rees Multiplicity Theorem for mixed multiplicities we need to consider a sequence of ideals.

\begin{Theorem} [Teissier, 1973] 
Let ${\bf I} = I_1, \ldots, I_s$ be a sequence of $\mm$-primary ideals.
For any $u = (u_1,...,u_s)\in \NN^s$ let ${\I}^u = I_1^{u_1}\cdots I_s^{u_s}$.
Then the  function $\ell(A/{\I}^u)$
is given by a  polynomial $P(u)$ of total degree $d$ 
for all $u \gg 0.$ 
The polynomial $P(u)$ can be  written  in terms of binomial coefficients as 
$$P(u)=\sum_{\a \in \NN^s,|\a| \leq d}e_\a({\I})
\binom{u_1+\a_1}{\a_1} \binom{u_2+\a_2}{\a_2}\cdots \binom{u_s+\a_s}{\a_s}
$$ 
where $e_\a({\I})$ are integers which  are positive if $|\a|=d.$
\end{Theorem}

The integers $e_\a(\I)$ with $|\a|=d$ are called the {\it mixed multiplicities} of the sequence  $\I.$
\cite{t1}.

Teissier showed that each mixed multiplicity of  $\I$ is  the multiplicity  of certain  ideals generated by  systems of parameters. 

Let $I = (c_1,...,c_r)$ be an arbitrary ideal in $A$. We say that a given property holds for a {\it sufficiently general element} $a \in I$ if there exists a non-empty Zariski-open subset $U \subseteq k^r$ such that whenever $a = \sum_{j=1}^r \alpha_jc_j$ and the image of $(\alpha_1,\ldots,\alpha_r)$ in $k^r$ belongs to $U$, then the given property holds for $a$.

\begin{Theorem}[Teissier, 1973] 
Let $\a = (\a_1,...,\a_s)$ with $|\a| = d$. Let $J$ be a parameter ideal generated by $\a_1$ general elements in $I_1$, ..., $\a_s$ general elements in $I_s$. Then
$$e_\a(\I)=e(J).$$
\end{Theorem}

This result was then generalized by Rees  for joint reductions \cite{r4}.  A sequence of elements $a_1, \ldots, a_s$ is called a {\it joint reduction} of the ideals $I_1,...,I_s$ in $A$
if $a_i \in I_i$ for $i=1,2,\ldots, s$ and the ideal $$\sum_{i=1}^{s}a_iI_1 \cdots I_{i-1}I_{i+1} \cdots I_s$$
is a reduction of $I_1 \cdots I_s.$
The existence of joint reductions was established first for a set of $d$ $\mm$-primary ideals by Rees and for 
any number of arbitrary ideals by O'Carroll \cite{oc}.

The next theorem is one of the fundamental results relating joint reductions with mixed multiplicities.

\begin{Theorem}[Rees Mixed Multiplicity Theorem, 1982] 
Let $\I$ be a sequence of $d$ $\mm$-primary ideals. Let $J$ be an ideal generated by a joint  reduction of  $\I$.
Put ${\bf 1}=(1,1,\ldots,1).$ Then
$$ e_{\bf 1}(\I)= e(J). $$
\end{Theorem}

It is now natural to ask for the converse of Rees Mixed Multiplicity theorem. This was done for  two-dimensional quasi-unmixed  local rings by Verma \cite{v1} and in any dimension by Swanson \cite{sw2}. 

\begin{Theorem}[Swanson's Mixed Multiplicity Theorem, 1992] 
Let $A$ be a quasi-unmixed local ring. 
Let $\I = I_1, \ldots, I_s$ be a sequence of ideals of $A.$ Let $a_j \in I_j$ for 
$j = 1, \ldots, s.$ Suppose the  ideals $(a_1, a_2, \ldots, a_s)$ 
and $I_1,  \ldots, I_s$ have equal radicals and their common height is $s.$ Suppose that
$$ e((a_1, \ldots, a_s)A_\wp) = e_{\bf 1} (I_1A_\wp,...,I_sA_\wp)
$$
for each minimal prime $\wp$ of $(a_1, \ldots, a_s).$ Then 
$a_1, \ldots, a_s$ is a joint reduction of $I_1, \ldots, I_s.$
\end{Theorem}

\section{Mixed multiplicities of two ideals: the general case} 

Now we are going to extend the notion of mixed multiplicities of two $\mm$-primary ideals $I,J$ to the case when
$I$ is an $\mm$-primary ideal and $J$ is an arbitrary ideal of a local ring $(A,\mathfrak m)$.
To this end  we associate with $I, J$ the standard bigraded algebra 
$$R(I|J) := \oplus_{u,v \ge 0} I^uJ^v/I^{u+1}J^v.$$
over the quotient ring $A/I$. 

Let $R = R(I|J)$. Since $A/I$ is an artinian ring,  $R$ has a  Hilbert polynomial $P_R(u,v)$.
We call the mixed multiplicities $e_{ij}(R)$ the {\it mixed multiplicities} of the ideals $I$ and $J$.

This notion coincides with the mixed multiplicities of the last section if $J$ is an $\mm$-primary ideal.
In fact, we have
$$\ell(A/I^uJ^v) = \sum_{t=0}^{u-1}\ell(I^tJ^v/I^{t+1}J^v)$$
for all $u, v \ge 0$. From this it follows that
$e_j(I|J) = e_{ij}(R)$ for $j < d$.

For this reason we will set
$$e_j(I|J) := e_{ij}(R)$$
for any $\mm$-primary ideal $I$ and any ideal $J$ in $A$.

Katz and Verma \cite{kv} showed if $\height J \ge 1$, then $\deg P_R(u,v) = \dim A -1$ and $e_0(I|J) = e(I)$.
This result can be generalized as follows \cite{Tr2}.

\begin{Lemma}\label{e0} 
Let $J$ be an  ideal. Then  $\deg P_R(u,v) = \dim A/0:J^\infty-1$ and   $e_0(I|J) = e(I,A/0:J^\infty)$. 
\end{Lemma}

We shall denote by $e(I,A/Q)$ the multiplicity of the ideal $(I+Q)/Q$ in the quotient ring $A/Q$ for any ideal $Q$ of $A$.

The positivity of the mixed multiplicities $e_j(I|J)$ is closely related to the dimension of the {\it fiber ring} of $I$, which is defined as the graded algebra
$$F(I) := \bigoplus_{n\ge 0}I^n/{\mathfrak m}I^n.$$
It can be shown that $J$ is a reduction of $I$ if and only if the ideal of $F(I)$ generated by the initial forms  of generators of $J$ in $F(I)_1 = I/\mm I$ is a primary ideal of the maximal graded ideal of $F(I)$. From this it follows that if the residue field of $A$ is infinite, the minimal number of generators of any minimal reduction $J$ of $I$ is equal to $\dim F(I)$. For this reason, $\dim F(I)$ is termed the {\it analytic spread} of $I$ and denoted by $s(I)$. If $I$ is an $\mm$-primary ideal, we have $s(I) = \dim A$. We refer the reader to \cite{nr} for more details. 

Katz and Verma \cite{kv} proved that $e_i(I|J) = 0$ for $i \ge s(J)$. This is a consequence of
the following bound for the partial degree $\deg_uP_R(u,v)$ of $P_R(u,v)$ \cite{Tr2}.  

\begin{Proposition} \label{deg u} 
$\deg_u P_R(u,v) < s(J)$. 
\end{Proposition}

\begin{Question}
Can one express $\deg_u P_R(u,v)$ in terms of $I$ and $J$?
\end{Question}  

To test the positivity of a mixed multiplicity $e_i(I|J) $ we have the following criterion, which also shows how to compute $e_i(I|J) $ as a Samuel multiplicity \cite{Tr2}.

\begin{Theorem}[Trung, 2001] \label{main} 
Let $J$ be an arbitrary ideal of $A$ and $0 \le i < s(J)$. Let $a_1,\ldots,a_i$ be elements in $J$ such that their images in $J/IJ$ and $J/J^2$ form filter-regular sequences in $R(I|J)$ and $R(J|I)$, respectively. Then $e_i(I|J) > 0$ if and only if  
$$\dim A/(a_1,\ldots,a_i):J^\infty = \dim A/0:J^\infty-i.$$ 
In this case, we have
$$e_i(I|J) = e(I,A/(a_1,\ldots,a_i):J^\infty).$$
 \end{Theorem}

Theorem \ref{main} requires the existence of elements in $J$ with special properties. However, if the residue field of $A$ is infinite, 
such elements always exist. In fact, any sequence of general elements $a_1,...,a_i$ in $J$ satisfies the assumption. 

As a consequence of Theorem \ref{main} we obtain the rigidity of mixed multiplicities and the independence of their positivity from the ideal $I$.

\begin{Corollary} \label{rigid3} 
Let $\rho = \max \{i|\ e_i(I|J) > 0\}.$ Then\par
{\rm (i) } $\height J -1 \le \rho \le s(J)-1$,\par
{\rm (ii) } $e_i(I|J) > 0$ for $0 \le i \le \rho$,\par
{\rm (iii)}  $\max \{i|\ e_i(I'|J) > 0\} = \rho$ for any $\mathfrak m$-primary ideal $I'$ of $A$. 
\end{Corollary}

Since $\rho$ doesn't depend on $I$, we  set $\rho(J) := \rho$.
In general we don't have the equality  $\rho(J) = s(J)-1$. 

\begin{Example} 
{\rm Let $A = k[[x_1,x_2,x_3,x_4]]/(x_1) \cap (x_2,x_3)$. Let $I$ be the maximal ideal of $A$ and $J = (x_1,x_4)A$. Then $F(J) \cong k[x_1,x_4]$. Hence $s(J) = \dim F(J) = 2$. One can verify that the images of $x_4$ in $J/IJ$ and $J/J^2$ are filter-regular elements in $R(I|J)$ and $R(J|I)$. We have $0:J^\infty = 0$ and $x_4A:J^\infty = (x_2,x_3,x_4)A$. Hence $\dim A/x_4A:J^\infty = 1 < \dim A/0:J^\infty - 1 = 2$. By Theorem \ref{main} this implies $e_1(I|J) = 0$.}
\end{Example}

We have the following sufficient condition for $\rho(J) = s(J)-1$. 

\begin{Corollary} \label{rigid4} 
Suppose all maximal chains of prime ideals in $A/0:J^\infty$ have the same length. Then $e_i(I|J)> 0$ for $0 \le i \le s(J)-1$. 
\end{Corollary}

\begin{Question}
Can one express $\rho(J)$ in terms of $J$?
\end{Question} 

For the computation of $e_i(I|J)$ we may replace $I$ and $J$ by their reductions. That means $e_i(I|J) = e_i(I'|J')$
for arbitrary reductions $I'$ and $J'$ of $I$ and $J$, respectively. 

Using reductions one obtains the following simple formula for the mixed multiplicities $e_i(I|J)$, $i \le \height J-1$ \cite{Tr2}.

\begin{Proposition} \label{parameter} 
Let $0 \le i \le \height J-1$. Let $a_1,\ldots,a_i$ and $b_1,\ldots,b_{d-i}$ be sufficiently general elements in $J$ and $I$, respectively. Then
$$e_i(I|J) =  e(I,A/(a_1,\ldots,a_i)) = e((a_1,\ldots,a_i,b_1,\ldots,b_{d-i})).$$ 
\end{Proposition}

As a consequence we obtain the following interpretation of $e_1(\mm|J)$ which was proved  by Katz and Verma \cite{kv1}.
For any ideal $J$ of $A$ we denote by $o(J)$ the $\mm$-{\it adic order} of $J$, that is, the largest integer $n$ such that $J \subseteq {\mathfrak m}^n$. 

\begin{Corollary} \label{e1} 
Let $(A,\mathfrak m)$ be a regular local ring and $J$ an ideal with $\height J \ge 2$. 
Then $$e_1({\mathfrak m}|J) = o(J).$$ 
\end{Corollary}

For $i \ge \height J$ we couldn't find a simple formula for $e_i(I|J)$ except in the following case \cite{Tr2}.

Recall that $J$ is called {\it generically a complete intersection} if $\height {\mathfrak p} = d-\dim A/J$ and $J_{\mathfrak p}$ is generated by $\height {\mathfrak p}$ elements for every associated prime ideal $\mathfrak p$ of $J$ with $\dim A/{\mathfrak p} = \dim A/J$. 

\begin{Proposition} \label{deviation} 
Let $J$ be an ideal of $A$ with $0 < s = \height J< s(J)$. Assume that $J$ is generically a complete intersection. Let $a_1,\ldots,a_s$ and $b_1,\ldots,b_{d-s}$ be sufficiently general elements in $J$ and $I$, respectively. Then
$$e_s(I|J) =  e(I,A/(a_1,\ldots,a_s)) - e(I,A/J).$$ 
\end{Proposition}

\section{Milnor numbers and mixed multiplicities} 

A geometric interpretation of the mixed multiplicities was found by  Teissier in the
 Carg\`{e}se paper \cite{t1} in 1973. Teissier was interested in Milnor numbers of isolated singularities of complex analytic hypersurfaces. We will now recall the concept of Milnor number and point out some of its basic properties found in \cite{mi}. 

Let $f : U \subset \CC^{n+1} \rightarrow  \CC$ be an analytic function in an open neighborhood $U$ of $\CC^{n+1}.$ Put 
$ S_\epsilon = \{z \in \CC^{n+1} : ||z|| = \epsilon\}.$  Define the map 
$\phi_\epsilon(z):S_\epsilon \setminus \{f=0\} \longrightarrow S^1$  by
$\phi_\epsilon(z)=f(z)/||f(z)||.$ Let $f_{z_i}:=\partial f/\partial z_i$ denote
the partial derivative of $f$ with respect to $z_i.$

Milnor and Palamodov proved the following

\begin{Theorem} [Milnor, Palamodov] 
If the origin is an isolated singularity of $f(z)$ then the fibers of 
$\phi_\epsilon$ for small $\epsilon$ have the  homotopy type of a 
bouquet of $\mu$ spheres of dimension $n$ having a single common point 
where 
$$\mu = \dim_{\CC} \frac{\CC\{z_0, z_1, \ldots , z_n\}}
{(f_{z_0}, \ldots, f_{z_n})}.$$
\end{Theorem}

The number $\mu$  is called the {\it Milnor number} of the isolated singularity. Therefore the Milnor number is nothing but the multiplicity of the Jacobian ideal 
$$J(f) := (f_{z_0}, \ldots, f_{z_n})$$
of $f.$  The Milnor number is a very useful invariant to detect the topology of a singularity.

Let $(X, x)$ and $(Y, y)$ be two germs of reduced complex analytic hypersurfaces of same dimension $n.$ 
Then they are called {\it topologically equivalent} if there exist representatives 
$(X, x) \subset  (U, x)$ and $(Y, y)\subset (V, y)$ where $U$ and $V$ 
are open in $\CC^{n+1}$ and a homeomorphism of pairs between $(U, x)$ and $(V, y)$ 
which carries $X$ to $Y.$

The basic result relating the Milnor number with the topology of the singularity is the following: 

\begin{Theorem}[Milnor, 1968] 
Let $(X, x)$ and $(Y, y)$ be two germs of hypersurfaces with isolated singularity having  same topological type. Then  $\mu_x(X)=\mu_y(Y).$
\end{Theorem}

Teissier \cite{t1} refined the notion of Milnor number.

\begin{Theorem}[Teissier, 1973] 
Let $(X, x)$ be a germ of a hypersurface in  $\CC^{n+1}$ with an isolated singularity. Let $E$ be an $i$-dimensional affine subspace of $\CC^{n+1}$ 
passing through $x.$ If $E$ is chosen sufficiently general then the Milnor number of $X \cap E$ at $x$ is independent of $E.$
\end{Theorem}

The Milnor number of $X \cap E,$ as in the above theorem, 
is denoted by $\mu^{(i)}(X, x).$ 

Note that $\mu^{(n+1)}(X,x)$ is the Milnor number of the isolated singularity. Moreover $\mu^{(1)}_x(X)=m_x(X)-1$ where $m_x(X)$ denotes the 
multiplicity of the hypersurface $X$ at $x.$ 
Put 
$$\mu^*(X, x) = (\mu^{(n+1)}(X, x), \mu^{(n)} (X, x), \ldots , \mu^{(0)}(X, x))$$
\medskip

\noindent
{\bf Teissier's Conjecture (cf. \cite{t1})} 
If $(X,x)$ and $(Y,y)$ have same topological type then
$$\mu^*_x(X)=\mu^*_y(Y).$$

The above conjecture contains Zariski's conjecture \cite{z} to the effect that $m_x(X)=m_y(Y)$ for topologically equivalent isolated singularities of hypersurfaces,  as a special case. Zariski's conjecture is still open. 

Suppose $f_t(z_0,z_1,\ldots, z_n)$ is an analytic
family of $n$-dimensional hypersurfaces with isolated singularity  at the origin. Suppose 
all these singularities have the same  Milnor number. Hironaka conjectured that under these 
hypotheses the singularities have same topological type when $n=1.$  This conjecture and the case of $n\not=2$ were  settled 
in the affirmative by Trang and  Ramanujam \cite{tr}.

A counterexample to Teissier's conjecture was given
in 1975 by J. Brian\c{c}on and J.-P. Speder by constructing a family of quasi-homogeneous
surfaces with contant Milnor number which is topologically trivial.

It is now known that the constancy of the sequence $\mu^*$ in a
family $F(t; x_1,\ldots ,x_n)=0$  of hypersurfaces with isolated
singularities at the origin is equivalent to the topological
triviality of general nonsingular sections of all dimensions through
the $t$-axis; this follows from \cite{t1}  and \cite{bs2}.

Teissier  devised a way to calculate the sequence $\mu^*(X, x).$ 
It turns out that the sequence $\mu^*(X, x)$  is identical to the sequence of mixed multiplicities of the Jacobian ideal $J(f )$ and $\mm$. 
More precisely the following result was proved by Teissier \cite{t1}.

Let $A=\CC\{z_0,z_1,\ldots,z_n\}$ denote the ring of convergent power series.
Let $f \in A$ be the equation of a hypersurface singularity $(X, 0).$ 
If $(X,0)$ is an isolated singularity, $J(f)$ is an $\mm$-primary ideal. Hence
the function $\ell(A/J(f )^u\mm^v)$ is given by a polynomial 
$P (u,v)$ of total degree $n + 1$ for large  $u$ and $v$. 

\begin{Theorem}[Teissier, 1973] 
Let $(X, 0)$ be a germ of a hypersurface in  $\CC^{n+1}$ with an isolated singularity. 
Put $\mu^{(i)}= \mu^{(i)}(X, 0).$  
Then the terms of total degree $n +1$ in $P(u,v)$ have the form
$$\frac{1} {(n+  1) \;!} \sum_{i=0}^{n+1}
\binom{n+1}{i}\mu^{(n+1-i)}u^{(n+1-i)}v^i.
$$
\end{Theorem}

If $(X,0)$ is not an isolated singularity, using mixed multiplicities in the general case we can also compute
the Milnor number of general hyperplane sections of $(X,0)$. 

Let $s$ be the codimension of the singular locus of $X$. Let $E$ be a general $i$-plane in $\CC^n$ passing through the origin, $i \le s$.
Then $(X \cap E,0)$ is an isolated singularity. Let $\mu(X \cap E,0)$ denote its Milnor number.

Let $a_1,\ldots,a_i$ be general elements in $J(f)$. Let $b_1,\ldots,b_{n-i}$ be the defining equation of $E$. 
It is easily seen that 
\begin{eqnarray*} 
\mu(X \cap E,0) & = &\ell(A/(a_1,\ldots,a_i,b_1,\ldots,b_{n-i}))\\
& = & e((a_1,\ldots,a_i,b_1,\ldots,b_{n-i})). 
\end{eqnarray*}
On the other hand, by Proposition \ref{parameter} we have 
$$e_i({\mathfrak m}|J(f)) = e((a_1,\ldots,a_i,b_1,\ldots,b_{n-i}))$$
for $i = 0,\ldots,s-1$. Hence we obtain the following formula for $\mu(X \cap E,0)$ \cite{Tr2}.

\begin{Theorem}
Let $(X,0)$  be a germ of a complex analytic hypersurface. With the above notations we have
$$\mu(X\cap E,0) = e_i({\mathfrak m}|J(f)).$$
\end{Theorem}

\section{Multiplicities of blow-up algebras} 

In this section we will discuss how mixed multiplicities arise naturally in the calculation of multiplicity of various blowup algebras. Let $(A,\mm)$ be a Noetherian local ring and $d = \dim A$. Let 
$I$ be an ideal of $A.$ The {\it Rees algebra} of $I$ is the graded $A$-algebra 
$$A[It] :=\bigoplus_{n=0}^{\infty} I^nt^n$$
where $t$ is an indeterminate. This graded algebra has a unique maximal homogeneous ideal
$M = (\mm, It).$ To compute the multiplicity of the local ring $A[It]_M$ we recall the following fact.

For any ideal $Q$ of a commutative ring $S$ we define the {\it associated graded ring} of $Q$ as the standard graded algebra
$$G_Q(S) = \bigoplus_{t\ge0}Q^n/Q^{n+1}$$
over the quotient ring $S/Q$. 
If $Q$ is a maximal ideal of $S$, then $G_Q(S)$ has the Hilbert function $\ell(Q^n/Q^{n+1})$.
It is easy to check that 
$$e(S_Q) = e(G_Q(S)).$$

Huneke and Sally \cite{hs} observed that
$$\ell(M^n/M^{n+1}) = \sum_{i=0}^{n} \ell (\mm^{n-i}I^i/\mm^{n-i+1}I^i).$$
They calculated this function for integrally closed $\mm$-primary ideals in a two dimensional regular local ring and obtained the formula
$$
e(A[It]_M ) = 1 + o(I),$$
where $o(I)$ denotes the $\mm$-adic order of $I$. 

On the other hand, we know by Corollary \ref{e1} that $e_1(\mm|I) = o(I)$. 
Therefore, the above formula can be rewritten as
$$e(A[It]_M ) = e_0(\mm|I)  + e_1(\mm|I).$$
That raises the natural question whether such a formula holds for an arbitrary ideal $I$. 

This question was answered by Verma in 1988  for $\mm$-primary ideals \cite{v2} and in 1992 for the general case \cite{v7}.

\begin{Theorem}[Verma, 1992] \label{Rees-Verma}
Let $(A, \mm)$ be a local ring of dimension $d.$ Let $I$ be an arbitrary ideal of positive height. Then
$$
e(A[It]_M ) = e_0(\mm|I) + e_1(\mm|I) +  e_2(\mm|I) + \cdots + e_{d-1}(\mm|I).$$ 
\end{Theorem}

As a consequence, the above result of Huneke and Sally also holds for an arbitrary $\mm$-primary ideal in a two dimensional regular local ring.

A similar formula for the multiplicity of the {\it extended Rees algebra}
$$A[It,t^{-1}] :=\bigoplus_{n \in \ZZ}^{\infty} I^nt^n\; (I^n = A\ \text{ for } n < 0)$$
was found by Katz and Verma \cite{kv}.  

\begin{Theorem}[Katz-Verma, 1992] \label{Katz-Verma}
Let $(A, \mm)$ be a local ring of dimension $d.$ Let $I$ be an ideal of positive height. 
Let $N = (\mm,It,t^{-1}) \subset A[It,t^{-1}]$.
Then
$$
e(A[It,t^{-1}]_N) = \frac{1}{2^d} \left[e(\mm^2 + I) + \sum^{d-1}_{j=0} 2^j e_j(\mm^2 + I|I)\right].
$$
\end{Theorem}

Similarly the multiplicity  of the fiber ring $F(I)$ can be expressed in terms  of a  mixed multiplicity of $\mm$ and $I$ \cite{jpv}. Furthermore, this knowledge also helps in detecting the Cohen-Macaulay property of $F(I)$ for certain classes  of ideals.

Though we have a well developed theory of mixed multiplicities when $I$ is an $\mathfrak m$-primary ideal \cite{t1}, \cite{r1}, there have been few cases where the mixed multiplicities can be computed in terms of well-known invariants of $\mathfrak m$ and $I$ when $I$ is not an $\mathfrak m$-primary ideal.

The first explicit computation of mixed multiplicities for a non-trivial case was done by Katz and Verma in \cite{kv1} for height 2 almost complete intersection prime ideals in a polynomial ring in three variables. 

In 1992, Herzog, Trung and Ulrich \cite{htu} computed the multiplicity of Rees algebras of homogeneous ideals generated by a $d$-sequence. 
Recall that a sequence of elements $x_1,\ldots,x_n$ is said to be a $d$-sequence if

(1) $x_i\notin (x_1,\ldots,x_{i-1},x_{i+1},\ldots,x_n),$

(2) $(x_1,\ldots,x_i):x_{i+1}x_k=(x_1,\ldots,x_i):x_k$ for all $k\ge i+1$ and all $i\ge 0$.

Examples of $d$-sequences are regular sequences, the maximal minors of an $n \times (n+1)$ matrix of indeterminates or the generators of an almost complete intersection \cite{Hu}. 

Let $A$ be a  standard graded ring over a field $k$  and $\mm$ the maximal graded ideal of $A$.
Let $I = (x_1, \ldots , x_n)$ be an ideal generated by a $d$-sequence of homogeneous elements in $R$ 
with $\deg (x_1)\leq \cdots \leq \deg(x_n).$  Herzog, Trung and Ulrich calculated the multiplicity of the Rees algebra $A[It]$ in terms of the multiplicities of $A/I_j$ where 
$$I_j = (x_1,\cdots, x_{j-1}) : x_j\;\; (j = 1, \ldots, n).$$ 
They used a technique which is similar to that of Gr\"obner bases and which does not involve mixed multiplicities.

\begin{Theorem}[Herzog-Trung-Ulrich, 1992] \label{HTR}
 Let $I$ be an ideal generated by a homogeneous $d$-sequence $x_1,\ldots,x_n$ of $A$ with $\deg x_1\le \ldots\le \deg x_n$. Then 
$$e(A[It]_M) = \left\{ \begin{array}{lll} 
\sum_{j=1}^se(A/I_j) & if & \dim A/I_1 = \dim A,\\
\sum_{j=1}^se(A/I_j) + e(A) & if & \dim A/I_1 = \dim A-1,\\
e(A) & if & \dim A/I_1 \le \dim A-2.
\end{array} \right.$$
\end{Theorem}

On the other hand, using essentially the same technique, N.D. Hoang \cite{Ho} was able to compute the mixed multiplicities
$e_i(\mm|I)$ in this case, namely, 
$$e_i(\mathfrak m|I)= \left\{ \begin{array}{ll} e(A/{I_{i+1}}) & \text{if $0 \le i\le s-1,$}\\
 0 & \text{if $i \ge s,$} 
\end{array}\right. $$ 
where  $s = \max \{i|\dim A/{I_i}=\dim A/{I_1}-i+1\}$. 
Combining this with Theorem \ref{Rees-Verma} he could recover the result of Herzog, Trung and Ulrich.

Using the same technique Raghavan and Verma \cite{rv} computed the Hilbert series of the associated graded ring 
$G_M(A[It]).$

\begin{Theorem}[Raghavan-Verma, 1997] 
Let $I$ be a graded ideal generated by a $d$-sequence as above.
Let $R =G_M(A[It])$. Then
$$
H(R; \lm_1, \lm_2) = H(A, \lm_1) + \lm_2
\sum_{s=1}^n \frac{H(A/I_s, \lm_1)} {(1 - \lm_2)^s}.
$$
\end{Theorem}

Inspired of work of Raghavan and Simis \cite{rasi}, Raghavan and Verma also computed the Hilbert series of the associated graded ring of homogeneous ideals generated by quadratic sequences, a generalization of $d$-sequences. As a consequence they obtained the following concrete formula for determinantal ideals.

\begin{Theorem}[Raghavan-Verma, 1997] 
Let $A$ be  the polynomial ring in $mn$ indeterminates over a field where $m \leq  n .$ Let $X$ be an $m \times  n$ 
matrix of these indeterminates. Let $I$ denote the ideal of 
$A$ generated by the maximal minors of $X.$ 
Then the Hilbert series of the  bigraded algebra $R =G_M(A[It])$ 
is given by the formula:
$$
H(R ; \lm_1, \lm_2) = H(A, \lm_1)  + \lm_2\sum_{\omega \in \Omega}
H (A/(\Pi^{\omega}), \lm_1) H(F_{\omega}, \lm_2)
$$
where 
\begin{eqnarray*}
\Pi  &=& \mbox{poset of all minors of all sizes of the matrix} \; X\\
\Omega &=& \mbox{ideal of}\; \; \Pi  \;\mbox{consisting of maximal minors of}\;\; X \\ \Pi^{\omega} &=& \{\pi \in \Pi : \pi \ngeq \omega\} \\ 
(\Pi^{\omega}) &=& \mbox{ideal of}\;  A \;\mbox{generated by} \;\;
\; \Pi^{\omega}.  \\
F_{\omega} &=& \mbox{the face ring over}\; k\; 
\mbox{of} \;\;
\Pi_{\omega} := \{\pi \in \Pi : \pi \leq \omega\}.
\end{eqnarray*}
\end{Theorem}

The techniques of Herzog, Trung and Ulrich was also used to compute the multiplicity of Rees algebras of ideals generated by filter-regular sequences of homogeneous elements \cite{Tr}. 

Recall that a  sequence $x_1,\ldots,x_n$ of elements of $A$ is called {\em filter-regular} with respect to $I$ if $x_i\notin \wp$ for all associated prime ideal $\wp \not\supseteq I$ of $(x_1,\ldots,x_{i-1})$, $i=1,\ldots,n$.
 
If $A$ is a {\it generalized Cohen-Macaulay ring}, that is, $A_\pp$ is Cohen-Macaulay and $\dim A/\pp + \height \pp = \dim A$ for all prime ideals $\pp \neq \mm$, then every ideal generated by a subsystem of parameters $x_1,\ldots,x_n$ is filter-regular with respect to  $I = (x_1,...,x_n).$

\begin{Theorem}[Trung, 1993] \label{filter-regular}
Let $I$ be a homogeneous ideal of $A$ generated by a subsystems of parameters $x_1,\ldots,x_n$  which is a filter-regular sequence with respect to  $I$ with $\deg x_1=a_1\le \ldots \le \deg x_n=a_n$. Then 
$$e\left(A[It]_M\right)=\left(1+\sum_{i=1}^{n-1}a_1\ldots a_i\right)e(A),$$
$$e\left(A[It,t^{-1}]_N\right)=\left(1+\sum_{i=l}^{n-1} a_1\ldots a_i\right)e(A),$$
where $l$ is the largest integer for which $a_l=1$, $l=0$ and $a_1\ldots a_l=1$ if $a_i>1$ for all $i=1,\ldots,n$.
\end{Theorem}

Comparing the first formula with Theorem \ref{Rees-Verma} one might ask whether in the above case,
$$e_i(\mathfrak m|I)= \left\{ \begin{array}{ll} a_1\cdots a_ie(A) & \text{if $0 \le i\le n-1,$}\\
 0 & \text{if $i \ge n.$} 
\end{array}\right. $$  
This has been proved by N. D. Hoang in \cite{Ho}.

\begin{Question}
Can one drop the condition $a_1 \le \ldots \le a_n$ in Theorem \ref{filter-regular}?
\end{Question} 

Note that the proof for a positive answer to this question in \cite{Vi2} is not correct.

Finally we report a formula for the multiplicity of the local ring $A[It]_{\mm A[It]}$. This is useful in finding when certain symmetric algebras are Cohen-Macaulay with minimal multiplicity.

\begin{Proposition}[Yoshida, 1995] 
Let $(R, \mm)$ be a quasi-unmixed local ring. Let $I$ be an ideal of 
positive height in $R$ with $\mu(I) =  s(I) = n.$ Then
$$ e\left(A[It]_{\mm A[It]}\right) 
= e_{n-1}(\mm \mid  I).$$
\end{Proposition}

\section{Mixed  multiplicities of a sequence of  ideals: the general case}

For some applications we need to extend the notion of mixed multiplicities to a sequence of ideals which are not necessarily $\mm$-primary.
This can be done in the same manner as for mixed multiplicities of two ideals.

Let $(A,\mm)$ be a local ring (or a standard graded algebra over a 
field with maximal graded ideal $\mm$). Let $I$ be an $\mm$-primary 
ideal and $\J = J_1,\ldots,J_s$ a sequence of ideals of $A$.
One can define the $\NN^{s+1}$-graded algebra
$$R(I|\J) := \bigoplus_{(u_0,u_1,...,u_s) \in 
\NN^{s+1}}I^{u_0}J_1^{u_1}...J_s^{u_s}/I^{u_0+1}J_1^{u_1}...J_s^{u_s}.$$
This algebra can be viewed as the associated graded ring of the ideal $(I)$ of
the {\it multi-Rees algebra} 
$$A[It_0,J_1t_1,...,J_st_s] := \bigoplus_{(u_0,u_1,...,u_s) \in 
\NN^{s+1}} I^{u_0}J_1^{u_1}...J_s^{u_s}t_0^{u_0}t_1^{u_1}...t_s^{u_s}.$$ 

For short, set $R = R(I|\J)$. Then $R$ is a standard 
$\NN^{s+1}$-graded algebra. Hence 
it has a Hilbert polynomial $P_R(u)$. For any $\a \in \NN^{s+1}$ with 
$|\a| = \deg P_R(u)$ we will set 
$$e_\a(I|\J) := e_\a(R).$$

If $J_1,\ldots,J_s$ are $\mm$-primary ideals, $e_\a(I|\J)$ 
coincides with the mixed multiplicities defined by using the function $\ell(A/I^{u_0}J_1^{u_1}...J_s^{u_s})$.
However, the techniques used in the $\mm$-primary 
case are not applicable for non $\mm$-primary ideals. For instance, mixed multiplicities of $\mm$-primary ideals are always positive, whereas they may be zero in the general case. 
We have to develop new techniques to prove the following general result which allows us to test the positivity of mixed multiplicities and  to compute them by means of Samuel's multiplicity \cite{TV}.

Throughout this section let  $J := J_1...J_s$ and  $d := \dim A/0:J^\infty.$

\begin{Theorem}[Trung-Verma, 2007] \label{dimension}
Let  $d \ge 1$. Then  $\deg P_R(u) = d - 1$ and
$$e_{(d-1,0,...,0)}(I|\J ) = e(I,A/0:J^\infty).$$
\end{Theorem}

We shall need the following notation for the computation of  mixed multiplicities.
A sequence of homogeneous elements $z_1,...,z_m$ in a multigraded algebra $S$ is 
called {\it filter-regular} if
$$[(z_1,...,z_{i-1}):z_i]_u = (z_1,...,z_{i-1})_u$$
for $u \gg 0$, $i = 1,...,m$. It is easy to see that this is equivalent to 
the condition
$z_i \not\in P$ for any associated prime 
$P \not\supseteq S_+$ of $S/(z_1,...,z_{i-1})$. 

We will work now in the $\ZZ^{s+1}$-graded algebra
$$S := \bigoplus_{u\in\ZZ^{s+1}}I^{u_0}J_1^{u_1}...J_s^{u_s}
/I^{u_0+1}J_1^{u_1+1}...J_s^{u_s+1},$$
which is the associated graded ring of 
the algebra $A[It_0,J_1t_1,...,J_st_s]$ with respect to the ideal $(IJ)$.

Let $\e_1,...,\e_m$  be any non-decreasing sequence of indices with 
$1 \le \e_i \le s$. Let $x_1,...,x_m$ be a sequence of elements of $A$ 
with $x_i \in J_{\e_i}$, $i = 1,...,m$. We denote by $x_i^*$  
the residue class of $x_i$ 
in $J_{\e_i}/IJJ_{\e_i}$.
We call $x_1,...,x_m$ an $(\e_1,...,\e_m)$-{\it superficial sequence} for the 
ideals $J_1,...,J_s$ (with respect to $I$) if $x_1^*,...,x_m^*$ is a 
filter-regular sequence in $S$. \par

The above notion can be considered as a generalization of  the classical notion of a 
superficial element of an ideal, which plays an important role in the theory 
of multiplicity. Recall that an element $x \in \aa $ is called superficial with respect to an 
ideal $\aa$ if there is an integer $c$ such that 
$$(\aa^n:x) \cap \aa^c = \aa^{n-1}$$ for 
$n \gg 0$.  A sequence of elements $x_1,...,x_m \in \aa$ is called a 
superficial sequence of $\aa$ if the residue class of $x_i$ in 
$A/(x_1,...,x_{i-1})$ is a superficial element of the ideal 
$\aa/(x_1,...,x_{i-1})$, $i = 1,...,m$. It is known that this is equivalent to 
the condition that the initial forms of $x_1,...,x_m$ in $\aa/\aa^2$ 
form a filter-regular sequence in the associated graded ring 
$\oplus_{n\ge 0}\aa^n/\aa^{n+1}$ (see e.g. \cite{Tr}).
\par

We have the following criterion for the positivity of mixed multiplicities
(a somewhat weaker result was obtained by Viet in \cite{Vi}). 

\begin{Theorem}[Trung-Verma, 2007] \label{positivity}
Let $\a = (\a_0,\a_1,...,\a_s)$ be any sequence of non-negative integers 
with $|\a| = d-1$.  Let $Q$ be any ideal generated by an 
$(\a_1,...,\a_s)$-superficial sequence of the ideals 
$I,J_1,...,J_s$.  Then $e_\a(I|\J) > 0$ 
if and only if   $\dim A/Q:J^\infty = \a_0+1.$ In this case, 
$$e_\a(I|\J) = e(I,A/Q:J^\infty).$$
\end{Theorem}

Let $k$ be the residue field of $A$. Using the prime avoidance 
characterization of a superficial element we can easily see that 
superficial sequences exist if $k$ is infinite. In fact, any sequence which consists  
of $\a_1$ general elements in $J_1$, ... , $\a_s$ elements in $J_s$ 
forms an $(\a_1,...,\a_s)$-superficial sequence  for the ideals $J_1,...,J_s$.

The following result shows that the positivity of mixed multiplicities does not depend on the ideal $I$ and that
the sequence of positive mixed multiplicities is rigid. 

\begin{Corollary} \label{rigid} 
Let $\a = (\a_0,\a_1,...,\a_s)$ be any sequence of non-negative integers 
with $|\a| = d -1$. Assume that $e_\a(I|\J) > 0$. 
Then\par
{\rm (a) } $e_\a(I'|\J) > 0$ for any $\mm$-primary ideal 
$I'$,\par
{\rm (b) } $e_\b(I|\J) > 0$ for all 
$\b = (\b_0,\ldots,\b_n)$ with $|\b| = d-1$ 
and $\b_i \le \a_i$, $i = 1,\ldots,n$.
\end{Corollary}

Mixed multiplicities of a sequence of ideals can be used to compute the multiplicity of multi-Rees algebras.

\begin{Theorem}[Verma, 1992] 
Let $\J = J_1, \ldots, J_s$ be a sequence of ideals of positive height. 
Let $M = (\mm, J_1t_1, \ldots, J_st_s) \subset A[J_1t_1, \ldots, J_st_s]$ and $e_1 = (1,0,...,0) \in \NN^{s+1}$.
Then
$$
e(A[J_1t_1, \ldots, J_st_s]_M) = 
\sum_{{\a}\in \NN^{s+1},\;\; |\a|=d-1} 
e_{\a + e_1}(\mm| \J). 
$$
\end{Theorem}

We shall see in the next section that mixed volumes of lattice polytopes are  special cases of mixed multiplicities of ideals.

\section{Mixed volume of lattice polytopes}

Let  us first recall the definition of mixed volumes.
Given two polytopes  $P, Q$ in $\RR^n$ (which need not to be 
different), their Minkowski sum is defined as the polytope
$$P + Q :=\{a + b \mid \ a \in P,\ b \in Q\}$$ 
The $n$-dimensional {\it mixed volume} of a collection of $n$ 
polytopes $Q_1,...,Q_n$ in 
$\RR^n$ is the value 
$$MV_n(Q_1, \ldots,Q_n) := 
\sum_{h=1}^s \sum_{1\le i_1 <...< i_h\le n} (-1)^{n-h}
  V_n(Q_{i_1}+\cdots + Q_{i_h}),$$
where $V_n$ denotes the $n$-dimensional Euclidean volume.
Mixed volumes play an important role in convex geometry \cite{BF} and elimination theory \cite{GKZ}.

Our interest in mixed volumes arises from the following result of Bernstein 
\cite{Be}, \cite{Kh}  which relates the number of solutions of a system of polynomial 
equations to the mixed volume of their Newton polytopes. 

For a Laurent polynomial $f \in \CC[x_1^{\pm 1},...,x_n^{\pm 1}]$ we denote by $Q_f$ the convex polytope spanned by the lattice points $\a = (\a_1,...,\a_n)$
such that the monomial $x_1^{\a_1}\cdots x_n^{\a_n}$ appears in $f$. This polytope is called the {\it Newton polytope} of $f$.

\begin{Theorem}[Bernstein, 1975] 
Let $f_1,...,f_n$ be Laurent polynomials in 
$\CC[x_1^{\pm 1},...,x_n^{\pm 1}]$  with 
finitely many common zeros in the torus $(\CC^*)^n,$  where $\CC^* = \CC \setminus \{0\}$.
Then the number of  common 
zeros of $f_1,...,f_n$ in $(\CC^*)^n$ is bounded above by the mixed volume
$MV_n(Q_{f_1},...,Q_{f_n})$. 
Moreover, this bound is attained  for a generic choice of coefficients in 
$f_1,..., f_n$.  
\end{Theorem}

Here, a generic choice of coefficients in  $f_1,..., f_n$ means that the 
supporting monomials of $f_1,..., f_n$ remain the same while their 
coefficients vary in a non-empty open parameter space. \par

Bernstein's theorem is a generalization of Bezout's theorem which says that if $f_1,...,f_n$ are polynomials  in $n$ variables having finitely many common zeros, then the number of common zeros of $f_1,...,f_n$ is bounded by $\deg f_1 \cdots\deg f_n$. In fact, by translation we may assume that the common zeros of $f_1,...,f_n$ lie in $(\CC^*)^n$. Let $P_i$ denote the $n$-simplex spanned by the origin and all points of the form $(0,...,\deg f_i,...,0)$. Then $Q_{f_i} \subseteq P_i$.
This implies $MV_n(Q_{f_1},...,Q_{f_n}) \le MV_n(P_1,...,P_n)$. It is easy to check that   $MV_n(P_1,...,P_n) = \deg f_1 \cdots\deg f_n$. 

Bernstein's theorem is a beautiful example of the interaction between  
algebra and combinatorics. The original proof in \cite{Be} has more or 
less a combinatorial flavor. A geometric proof using intersection theory was 
given by Teissier \cite{t5} (see also the expositions \cite{Fu}). Here we  sketch 
an algebraic proof of Bernstein's theorem by means of mixed multiplicities.

First of all, using homogenization we can reformulate Bernstein's theorem as follows.

\begin{Theorem} \label{homogen} 
Let $k$ be an algebraically closed field. Let $g_1,..., g_n$ be homogeneous 
Laurent polynomials in $\CC[x_0^{\pm 1},x_1^{\pm 1},...,x_n^{\pm 1}]$  with 
finitely many common zeros in $\PP_{\CC^*}^n$.  Then 
$$ | \{ \a \in \PP_{\CC^*}^n \mid g_i(\a)=0,\; i=1, 2, \ldots,n \}|
   \leq \frac{MV_n(Q_{g_1},...,Q_{g_n})}{\sqrt{n+1}}.$$
Moreover, this bound is attained for a generic choice of coefficients in  
$g_1,..., g_n$.
\end{Theorem}

We may reduce the above theorem to the case of polynomials. In fact, 
if we multiply the given Laurent polynomials with an appropriate  
monomial, then we will obtain a new system of polynomials. Obviously, the new 
polynomials  have the same common 
zeros in $\PP_{\CC^*}^n$. Since their Newton polytopes are translations of the old ones, their 
mixed volumes do not change, too.

Now assume that $g_1,...,g_n$ are homogeneous polynomials in the polynomial ring $A  = 
k[x_0,x_1,...,x_n]$, where $k$ is a field. 
Let $M_i$ be the set of monomials occuring in $g_i$.
Let $\mm$ be the maximal graded ideal of $A$ and $J_i$ the 
ideals of $A$ generated by $M_i$. 
Put 
$$R = R(\mm|J_1,...,J_n).$$
We know by Theorem \ref{dimension} that $\deg P_R(u) = n$.

First, using the interpretation of mixed multiplicities as Samuel's multiplicity we 
can prove the following bound for the number of common zeros of $g_1,...,g_n$ 
in $\PP_{k^*}^n$, where $k^* = k \setminus \{0\}$.

\begin{Theorem}[Trung-Verma, 2007] \label{multiplicity} 
Let $k$ be an algebraically closed field. Let $g_1,..., g_n$ be homogeneous 
polynomials in $k[x_0,x_1,...,x_n]$  with 
finitely many common zeros in $\PP_{k^*}^n$. Then 
$$ 
| \{ \a \in \PP_{k^*}^n \mid g_i(\a)=0,\; i=1, 2, \ldots,n \}|
 \leq  e_{(0,1,...,1)}(R). 
$$
Moreover, this bound is attained for a 
generic choice of coefficients in  $g_1,..., g_n$ if $k$ has 
characteristic zero.
\end{Theorem}

It remains to show that 
$$e_{(0,1,...,1)}(R) = \frac{MV_n(Q_{g_1},...,Q_{g_n})}{\sqrt{n+1}}.$$

To prove that we shall need the following basic property of mixed volumes.

Let $g_0 = x_0  \cdots  x_n$ and $\Q = (Q_{g_0},Q_{g_1},...,Q_{g_n})$.
Let $\l = (\l_0,...,\l_n)$ be any sequence of positive integers.
We denote by $\l\Q$ the Minkowski sum $\l_0Q_0+  \cdots + \l_nQ_n$ 
and by
$\Q_\l$ the multiset of $\l_0$ copies of polytopes $Q_0$,...,$\l_n$ copies of polytopes 
$Q_n$.
Minkowski showed that the volume of the polytope $\l\Q$ is a 
homogeneous polynomial in $\l$ whose coefficients are mixed volumes up to 
constants:
$$V_n(\l\Q) = \sum_{\a \in \NN^s, |\a| = 
n}\frac{1}{\a!}MV_n(\Q_\a)\l^\a.$$

On the other hand, there is also a Minkowski formula for mixed multiplicities,
which  arises  in the computation of the multiplicity of the $\NN$-graded algebra
$$R^\l := \bigoplus_{t \ge 0}R_{t\l}.$$
One calls $R^\l$ the $\l$-{\it diagonal subalgebra} of $R$. This notion plays 
an important role in the study of embeddings of blowups of projective 
schemes \cite{CHTV}. \par

It is easy to check that for all positive integers $\l,$
$$e(R^\l) = n!\displaystyle \sum_{\a \in \NN^{s+1},\;\; |\a| = n} 
\frac{1}{\a!}e_\a(R)\l^\a.$$

Since $R^\l$ is a ring generated by monomials of the same degree, using Ehrhart's theory for the  
number of lattice points in  lattice polytopes 
(see  e.g. \cite{St}) we have
$$e(R^\l) = \frac{n!V_n(\l\Q)}{\sqrt{n+1}}.$$

Now we can compare the two Minkowski formulas and obtain the following relationship between mixed multiplicities and mixed volumes.

\begin{Theorem}[Trung-Verma, 2007]
With the above notation we have 
$$e_\a(R) = \frac{MV_n(\Q_\a)}{\sqrt{n+1}}$$
for any $\a \in \NN^{n+1}$ with $|\a| = n$.
\end{Theorem}

As a consequence, for $\a = (0,1,...,1)$ we get 
$$e_{(0,1,...,1)}(R) = \frac{MV_n(Q_{g_1},...,Q_{g_n})}{\sqrt{n+1}}= MV_n(Q_{f_1},...,Q_{f_n}),$$
which completes the proof of Bernstein's theorem.  Similarly, we can show that Bernstein's theorem holds for polynomials over any algebraically closed field of characteristic zero.

Since any collection of $n$ lattice polytopes in $\RR^n$ can be realized as the Newton polytopes of $n$ polynomials,
we can always express mixed volumes of lattice polytopes as mixed multiplcities of ideals.

It is known that computing  mixed volumes is a hard enumerative problem. 
Instead of that we can now compute mixed multiplicities of 
the associated graded ring of the multigraded Rees algebra 
$A[J_1t_1,...,J_nt_n]$ with respect to the ideal $\mm$. By 
Theorem \ref{positivity}, these mixed multiplicities can be interpreted as 
Samuel multiplicities. The computation of these multiplicities can be 
carried out by computer algebra systems such as {\it Cocoa, Macaulay 2} and 
{\it Singular.}

\section{Minkowski inequalities and equalities}

Let $(A,\mm)$ be a local ring of dimension $d$. Suppose $I$ and $J$  are $\mm$-primary ideals. Then $e(IJ)$ can be computed 
if the mixed multiplicities of $I$ and $J$ are known:
$$
e(IJ) = e(I)  +\binom{d}{1} e_1(I|J)  +\cdots + \binom{d}{i}e_i(I|J) + \cdots + e(J).$$ 

Teissier \cite{t1} made the following conjectures  based on the comparison of this formula with the binomial expansion
$$(e(I)^{\frac{1}{d}} + e(J)^{\frac{1}{d}})^d=
e(I)  +\cdots + \binom{d}{i}e(I)^{\frac{d-i}{d}}e(J)^{\frac{i}{d}}+\cdots +e(J).$$
\sk

\noindent {\bf Teissier's First Conjecture}:
{\em For all $i = 0, 1, ...,  d,$}
$$e_i(I|J)^d \leq e(I)^{d-i}e(J)^i.$$

The validity of Teissier's First Conjecture implies the 
{\em Minkowski's inequality}  for multiplicities:
$$e(IJ)^{1/d} \leq e(I)^{1/d} + e(J)^{1/d}.$$ 
\sk

\noindent {\bf Teissier's Second Conjecture}: 
{\em Let $d \geq 2.$ Put $e_i(I|J)=e_i$ for $i=1, \ldots, d-1.$ Then}
$$\frac{e_1}{e_0} \leq \frac{e_2}{e_1} \leq \cdots \leq \frac{e_d}{e_{d-1}}.$$

Teissier proved that the second conjecture implies the first. He also showed the validity of the second conjecture for reduced Cohen-Macaulay complex analytic algebras  \cite{t3}. Rees and Sharp \cite{resh} investigated these conjectures for all local rings
and proved:

\begin{Theorem}[Rees and Sharp, 1978] \label{Rees-Sharp}
Teissier's conjectures  and hence Minkowski inequality for multiplicities holds  for all local rings.
\end{Theorem}

It is natural to ask when equalities hold in Minkowski inequalities. It is easy to see that if $I$ and $J$ are {\it projectively equivalent}, that is, there exist positive integers $r$ and $s$ so that $I^r$ and $J^s$ have the same integral closures, then 
$$e(IJ)^{1/d} = e(I)^{1/d}+  e(J)^{1/d}.$$
The converse was proved by Teissier \cite{t4} for Cohen-Macaulay normal complex analytic algebras by using mixed multiplicities. Then Katz \cite{ka} showed that  in quasi-unmixed local rings,
Minkowski equalities hold for $\mm$-primary ideals $I$ and $J$ if and only if 
they are projectively equivalent. 

It is interesting to note that  Rees Multiplicity Theorem is a consequence of Minkowski equalities. 
We reproduce Teissier's lightning proof  of the converse of Rees Multiplicity Theorem found in \cite{t4}.

Let $A$ be a quasi-unmixed local ring. Let $J \subseteq I$ be $\mm$-primary ideals with $e(I)=e(J)=e.$ Since $JI \subseteq I^2,$ $e(IJ) \geq e(I^2)=2^d e(I).$ Hence
$$2 e^{1/d} \leq e(IJ)^{1/d} \leq e^{1/d}+  e^{1/d} =2 e^{1/d},$$
which implies $e(IJ)^{1/d} = e^{1/d}+  e^{1/d}$. Therefore, $I$ and $J$ are projectively equivalent
so that there exist positive integers $r$ and $s$ such that $\overline{I^r}=\overline{J^s}.$
It follows that $e(I^r) = e(J^s)$. Since $e(I^r)=r^de$ and $e(J^s)=s^de,$ we get $r = s$ which means that $J$ is a reduction of $I.$

We have seen in the preceding section that mixed volume is only a special case of mixed multiplicities.
Threfore, properties of mixed volumes may predict unknown properties of mixed multiplicities. For instance, 
consider the famous Alexandroff-Fenchel inequality among mixed 
volumes:
$$MV_n(Q_1,...,Q_n)^2 \ge MV_n(Q_1,Q_1,Q_3,...,Q_n)MV_n(Q_2,Q_2,Q_3,...,Q_n).$$
Khovanski \cite{Kh} and Teissier \cite{t5} used the Hodge index theorem in 
intersection theory to prove this inequality.
This leads us to believe that a similar inequality should hold for 
mixed multiplicities \cite{TV}. 

\begin{Question} 
{\rm Let $(A,\mm)$ be a local (or standard graded) ring 
with $\dim A = n+1 \ge 3$. Let $I$ be an $\mm$-primary ideal and 
$J_1,...,J_n$ ideals of height $n$. Put $\alpha=(0,1, \ldots, 1).$
Is it true that
$$
e_{\alpha}(I|J_1,...,J_n)^2 \ge  
e_{\alpha}(I|J_1,J_1,J_3,...,J_n)e_{\alpha}(I|J_2,J_2,J_3,...,J_n)~?
$$}
\end{Question}

Using Theorem \ref{positivity} we can reduce this theorem to the case 
$\dim A = 3$.
In this case, we have to prove the simpler formula:
$$e_{(0,1,1)}(I|J_1,J_2)^2 \ge e_{(0,1,1)}(I|J_1,J_1)e_{(0,1,1)}(I|J_2,J_2).$$
Unfortunately, we were unable to give an answer to the above question.

The difficulty can be seen from the fact that 
the above inequality does not hold if $J_1,...,J_n$ are 
$\mm$-primary ideals. In fact,  we can even show that the 
reverse  inequality holds, namely,
$$e_{\alpha}(I|J_1,...,J_n)^2 \le  
e_{\alpha}(I|J_1,J_1,J_3,...,J_n)e_{\alpha}(I|J_2,J_2,J_3,...,J_n)$$
where $\alpha=(0,1,1, \ldots, 1).$
For that we only need to show the inequality
$$e_{(1,1)}(J_1|J_2)^2 \le e(J_1,A)e(J_2,A),$$
for a two-dimensional ring $A$. But this is a special case of Theorem \ref{Rees-Sharp}.

\section{The multiplicity sequence}

The main aim of this section is to present another generalization of the multiplicity of an ideal by mixed multiplicities.

Let $(A,\mm)$ be a local ring of dimension $d$ and $I$ an ideal of $A$.
If $I$ is $\mm$-primary, we can consider the Hilbert-Samuel function $\ell(A/I^t)$ and define the multiplicity $e(I)$. Actually, $e(I)$ is the multiplicity of the associated graded ring $G := G_I(A)$.

If $I$ is not $\mm$-primary, we can replace $G$ by the associated graded ring of the ideal $\mm G$ of $G$:
$$R := G_{\mm G}(G) = \bigoplus_{(u,v) \in \NN^2}\big(\mm^uI^v + I^{v+1}/\mm^{u+1}I^v +I^{v+1}\big).$$
This is a standard bigraded algebra over the residue field $A/\mm$. 
Hence we can consider the Hilbert function $H_R(u,v)$. 
Since $H_R(u,v)$ is a polynomial for $u,v$ large enough, the sum transform
$$H_R^{(1,1)}(u,v) := \sum_{i = 0}^u\sum_{j = 0}^vH_R(u,v)$$
is given by a polynomial $P_R^{(1,1)}(u,v)$ for $u,v$ large enough. It is easy to check that $\deg P_R^{(1,1)}(u,v)  = d$. 
If we write this polynomial in the form
$$P_R^{(1,1)}(u,v) =\sum_{i= 0}^d 
\frac{c_{i \; d-i }(R)}{i!(d-i)!}u^iv^{d-i}  +
\text{\rm lower-degree terms},$$
then $c_{i\;d-i}(R) $ are non-negative integers for $i = 0,...,d$. We set 
$$c_i(I) := c_{i\;d-i}(R).$$

It is easily seen that the mixed multiplicities of $R$ belong to the multiplicity sequence: $e_{ij}(R) = c_{i+1}(I)$ for $i+j = d-2$.

The multiplicity of the associated graded ring $G$ with respect to the maximal graded ideal can be expressed as the sum of the multiplicity sequence.

\begin{Theorem}[Dade, 1960]
Let $M$ denote the maximal graded ideal of the associated graded ring $G$ of $I$. Then
$$e(G_M) = \sum_{j=0}^dc_j(I).$$
\end{Theorem}

Achilles and Manaresi \cite{AM2} call $c_0(I),...,c_d(I)$ the {\it multiplicity sequence} of $I$.
The multiplicity sequence can be considered as a generalization of the multiplicity of an $\mm$-primary ideal. 

\begin{Theorem}[Achilles-Manaresi, 1997]
Let $s = s(I)$ and $r = \dim A/I$. Then\par
{\rm (a) } $c_j(I) = 0$ for $j <d-s$ and $j > r$,\par
{\rm (b) } $c_{d-s}(I) = \sum e(\mm G_P)e(G/P),$ where $P$ runs through all highest associated prime ideals of $\mm G$ such that $\dim G/P + \height P = d$, \par
{\rm (c) } $c_r(I) = \sum e(I_\wp)e(A/\wp),$ where $\wp$ runs through all highest associated prime ideals of $I$ such that $\dim A/\wp + \height \wp = d$.
\end{Theorem}

As a consequence, if $I$ is an $\mm$-primary ideal, then $c_0(I) = e(I)$ and $c_i(I) = 0$ for $i > 0$.

In particular, $c_0(I) > 0$ if and only if $s(I) = d$. In this case, we set $j(I) := c_0(I)$ and call it the $j$-{\it multiplicity} of $I$ \cite{AM1}. 

Using $j$-multiplicity one can extend Rees Multiplicity Theorem for arbitrary ideals as follows \cite{FM}.

\begin{Theorem}[Flenner-Manaresi-Ulrich]
Let $J$ be an ideal in $J$.
If $J$ is a reduction of $I$, then $j(J_\wp) = j(I_\wp)$ for all prime ideals $\wp \supseteq I$ with $s(I_\wp) = \dim A_\wp$. The converse holds if $A$ is a quasi-unmixed local ring.
\end{Theorem}

The multiplicity sequence can be computed by the following formula.

\begin{Theorem}[Achilles-Manaresi, 1997] \label{AM}
 Let $Q = (x_1,...,x_s)$ be a minimal reduction of $I$ such that the images of $x_1,...,x_s$ in $R_{(0,1)} = I/\mm I$ is a filter-regular sequence of $R$ with respect to the ideal $(R_{(0,1)})$.
Set $Q_0 := 0$ and $Q_i = (x_1,...,x_i)$ for $i = 1,...,s$. Then
$$c_{d-i}(I) = \sum \ell(A_\wp/(Q_{i-1}:I^\infty,x_i)_\wp)e(A/\wp),$$
where $\wp$ runs through all associated prime ideals containing $I$ of the ideal $(Q_{i-1}:I^\infty,x_i)$.
\end{Theorem}

For the computation of the multiplicity sequence we may assume that the residue field of $A$ is infinite.
In this case we can always find a minimal reduction $Q$ of $I$ which satisfies the assumption of the above theorem.
Moreover, if $J$ is a reduction of $I$, then ideal $J^*$ of $R$ generated by the images of the elements of $J$  in $R_{(0,1)}= I/\mm I$ has the property $(J^*)_{(u,v)} = (R_{(0,1)})_{(u,v)}$ for all $u,v$ large enough. 
Therefore, we can find a minimal reduction $Q$ of $J$ (and hence of $I$) which satisfies the assumption of the above theorem for both ideals $J$ and $I$. As an immediate consequence we obtain the following fact 
(a complicated proof was given by Ciuperca in \cite{C}): 

\begin{Corollary}
Let $J$ be a reduction of $I$, then $c_i(J) = c_i(I)$ for all $i = 0,...,d$.
\end{Corollary}

Inspired of Rees Multiplicity Theorem we raise the following question.

\begin{Question}
{\rm Let $A$ be a quasi-unmixed local ring and $J$ an ideal in $I$ with $\sqrt{J} = \sqrt{I}$. 
Is $J$ a reduction of $I$ if $c_i(J) = c_i(I)$ for $i = 0,...,d$?}
\end{Question}

The multiplicity sequence can be used to compute  the degree of the St\"uckrad-Vogel cycles in the intersection algorithm.
 
Let $X, Y \subset {\PP}_k^n$ be two equidimensional subschemes. One can associate with $X \cap Y$ certain cycles $v_1,...,v_n$ as follows \cite{SV}, \cite{Vo}.
Let $V$ be the ruled join variety of $X$ and $Y$ in 
$${\PP}_{k(t)}^{2n+1} = \Proj\ k(t)[x_0,\ldots,x_n,y_0,\ldots,y_n].$$
where $k(t) = k(t_{ij}|\ 1 \le i \le n+1, 0 \le j \le n)$ is a pure transcendental extension of $k$. 
Put $w_0 = [V]$. Let $E$ be the linear subspace of 
${\PP}_{k(t)}^{2n+1}$ given by the equations $x_0-y_0 = \cdots = x_n-y_n = 0$. For $i = 0,\ldots,n-1$ let $h_i$ denote the divisor of $V$ given by the equation $\sum_{j=0}^nt_{ij}(x_j-y_j) = 0$. If $w_{i-1}$ is defined for some $i \ge 1$, we decompose $w_{i-1} \cap h_{i-1} = v_i+w_i,$ where the support of $v_i$ lies in $E$ and $w_i$ has no components contained in $E$. 
Using the cycles $v_1,...,v_n$ St\"uckrad and Vogel proved that there exist a set ${\Lambda (X,Y)}$ of irreducible subschemes $C$ of $(X \cap Y) \times_k k(t)$  and intersection numbers $j(X,Y;C)$ such that
$$\deg X\deg Y = \sum_{C \in {\Lambda(X,Y)}}j(X,Y;C)\deg C.$$
Algebraically, if we set $A = k(t)[x_0,\ldots,x_n,y_0,\ldots,y_n]/(I_X,I_Y)$, where $I_X$ and $I_Y$ denote the defining ideals of $X$ and $Y$ in $k[x_0,\ldots,x_n]$ and $k[y_0,\ldots,y_n]$,  and $I = (x_0-y_0,\ldots,x_n-y_n)A$, 
then $\deg v_i = c_{d-i}(I)$ by Theorem \ref{AM}. \par

Using Theorem \ref{main} we can also describe  $\deg v_i$ in terms of the mixed multiplicities $e_i({\mathfrak m}|I)$
\cite{Tr2}.

\begin{Theorem}[Trung, 2003]
With the above notations we have
$$\deg v_i = e_{i-1}({\mathfrak m}|I) - e_i({\mathfrak m}|I).$$ 
\end{Theorem}

Achilles and Rams \cite{AR} showed that the Segre numbers introduced by Gaffney and Gassler \cite{GG} in singularity theory and the extended index  of intersection introduced by Tworzewski \cite{Tw} in analytic intersection theory are special cases of the multiplicity sequence. We refer the readers to the report \cite{AM3} for further applications of the multiplicity sequence.

\section{Hilbert function of non-standard bigraded algebras}

In general, the Hilbert function $H_R(u,v)$ of a finitely generated bigraded algebra $R$ over a field $k$ is not a polynomial for large $u,v$.

In this section we will study the case when $R$ is generated by elements of 
bidegrees $(1,0),(d_1,1),\ldots,(d_r,1)$, 
where $d_1,\ldots,d_r$ are non-negative integers. 
This case was considered first by P.~Roberts in \cite{ro} where it is
shown that there exist integers $c$ and $v_0$ such that $H_R(u,v)$ 
is equal to a polynomial $P_R(u,v)$ for $u \ge cv$ and $v \ge v_0$. 
He calls $P_R(u,v)$ the {\it Hilbert polynomial} of the bigraded 
algebra $R$.

It is worth remarking  that Hilbert polynomials of bigraded algebras of 
the above type appear in Gabber's proof of Serre's non-negativity 
conjecture (see e.g. \cite{Ro2}) and that the positivity of certain coefficient 
of such a Hilbert polynomial is strongly related to Serre's positivity
conjecture on intersection multiplicities \cite{Ro3}. \par

Roberts' result can be made more precise as follows \cite{HT}.

\begin{Theorem} \label{exist} 
Let $d = \max\{d_1,\ldots,d_r\}.$
There exist integers $u_0,v_0$ such that for  $u\ge dv+ u_0$ and $v \ge 
v_0$,  $H_R(u,v)=P_R(u,v)$.
\end{Theorem}

If $R$ is a standard bigraded algebra, then $d = 0$. Hence the Hilbert function $H_R(u,v)$ is given by a polynomial for $u,v$ large enough. As in the standard bigraded case, the total degree $\deg P_R(u,v)$ can  also be expressed in terms of the relevant dimension of $R$ \cite{HT}.

\begin{Theorem}  \label{total 1}
$\deg P_R(u,v) = \rdim R-2.$
\end{Theorem}

The partial degree $\deg_u P_R(u,v)$ can be expressed in terms of the graded modules
$$R_v := \bigoplus_{u \ge 0}R_{(u,v)}.$$
Note that $R_0$ is a finitely generated standard $\NN$-graded algebra
and $R_v$ is a finitely generated graded $R_0$-module. Define
$$\sdim R := \dim (R/0:R_+^\infty)_0.$$
If $R$ is a standard bigraded algebra, then $(R/0:R_+^\infty)_0 = R/(0:R_+^\infty + (R_{(0,1)}))$.
 
\begin{Theorem} \label{partial} 
For $t$ large enough,
$$\deg_u P_R(u,v)  = \dim R_t-1 = \sdim R.$$
\end{Theorem} 

This result was already proved implicitly by P. Roberts for bigraded algebras generated by 
elements of bidegree $(1,0),(0,1),(1,1)$ \cite{Ro3}. 

By Theorem \ref{total 1} and Theorem \ref{partial} we always have
$$\rdim R = \deg P_R(u,v) \ge \deg_u P_R(u,v) = \sdim R.$$
Note that the inequality may be strict.

\begin{Question}
Does there exist similar formulas for $\deg_v P_R(u,v) $?
\end{Question}

Now we write the Hilbert polynomial $P_R(u,v)$ in the form
$$P_R(u,v) = \sum_{i= 0}^s \frac{e_i(R)}{i!(s-i)!}u^iv^{s-i}  +
\text{\rm lower-degree terms},$$
where $s = \deg P_R(u,v)$. Following Teissier we call the numbers $e_i(R)$ the {\it mixed multiplicities} of 
$R$. One can show that the mixed multiplicities $e_i(R)$ satisfy the associativity formula of Proposition \ref{associative}.

Unlike the case of standard bigraded algebras, a mixed multiplicities $e_i(R)$ may be negative.

\begin{Example} \label{polynomial} 
{\rm Let $S = 
k[X_1,\ldots,X_m,Y_1,\ldots,Y_n]$ $(m \ge 1, n \ge 1)$ be a bigraded polynomial ring with 
$\deg X_i = (1,0)$ and $\deg Y_j  = (d_i,1).$
We have $H_S(u,v) = P_S(u,v)$ for $u \ge dv$ with $\deg  P_S(u,v)  = m+n-2$ and
$$e_{i,m+n-2-i} = \left\{\begin{array}{ll} (-1)^{m-i-1}\displaystyle 
\sum_{j_1+ \ldots + j_n = m-1-i}d_1^{j_1}\cdots d_n^{j_n} & \text{if }\ i 
< m,\\
0 & \text{if }\ i \ge m.
\end{array}\right.$$}
\end{Example}

We set  $\rho_R := \max\{i|\ e_i(R) \neq 0\}.$

\begin{Theorem}[Hoang-Trung, 2003] \label{coefficient}
The mixed multiplicities $e_i(R)$ are integers with $e_{\rho_R}(R) > 
0$.
\end{Theorem}

We always have $\rho_R \le \deg_uP_R(u,v)$.  
The following result gives a sufficient condition for $\rho_R = \deg_uP_R(u,v)$. 
Note that this condition is satisfied if 
$R$ is a domain or a Cohen-Macaulay ring.

\begin{Proposition}[Hoang-Trung, 2003] \label{equi}
Suppose $\dim R/P = \rdim R$ for all minimal prime ideals of $\Proj 
R$. Let $d = \deg_u P_R(u,v)$. Then $e_d(R) > 0$.
\end{Proposition}

As a consequence we obtain the following generalization of a result  by P. Roberts \cite{Ro3} if 
$R$ is generated by elements of degree $(1,0),(0,1),(1,1)$.
That result was used to give a criterion for the positivity of Serre's 
intersection multiplicity.

\begin{Corollary}  
Suppose  there exists an associated prime ideal $P$ with $\dim R/P = \rdim R$
and $\sdim R/P = \sdim R$. Let $s = \sdim R$. Then $e_s(R) > 0$.
\end{Corollary}

\begin{Question}
Can one describe $\rho_R$ in terms of well-understood invariants of $R$?
\end{Question}

\section{Hilbert function of bigraded Rees algebras}

The inspiration for our study on Hilbert function of non-standard bigraded algebras
comes mainly from the fact these algebras include Rees algebras of homogeneous ideals.

Let $A$ be a standard graded algebra over a field $k$. Let $I$ be a 
homogeneous ideal of $A$. The Rees algebra  $A[It]$ is naturally bigraded:
$$A[It]_{(u,v)} := (I^v)_ut^v$$
for all $(u,v) \in \NN^2$. 

Let $A = k[x_1,\ldots,x_n]$, where 
$x_1,\ldots,x_n$ are homogeneous elements with $\deg x_i = 1$. Let $I = 
(f_1,\ldots,f_r)$, where $f_1,\ldots,f_r$ are homogeneous elements with 
$\deg f_j = d_j$. Put $y_j = f_jt$. Then $A[It]$ is generated by the 
elements $x_1,\ldots,x_n$ and $y_1,\ldots,y_r$ with $\deg x_i = (1,0)$ 
and $\deg y_j = (d_j,1)$. Hence $A[It]$ belongs to the class of 
bigraded algebras considered in the preceding section.

\begin{Theorem}[Hoang-Trung, 2003] \label{Rees} 
Set $d = \max\{d_1,\ldots,d_r\}$
and $s = \dim A/0:I^\infty-1$.  There exist  integers $u_0,v_0$ such 
that for  
$u \ge dv+u_0$ and $v \ge v_0$, the Hilbert function $H_{A[It]}(u,v)$ 
is equal to a polynomial $P_{A[It]}(u,v)$ with 
$$\deg P_{A[It]}(u,v) = \deg_u P_{A[It]}(u,v) = s.$$
Moreover, if $s \ge 0$ and $P_{A[It]}(u,v)$ is written in the form 
$$P_{A[It]}(u,v) = \sum_{i=0}^s\frac{e_i(A[It])}{i!(s-i)!}u^iv^{s-i}  +
\text{\rm lower-degree terms},$$
then the coefficients $e_i(A[It])$ are integers for all $i$ with 
$e_s(A[It])  = e(A/0:I^\infty)$. 
\end{Theorem}

This result has some interesting applications.
First of all, the Hilbert polynomial $P_{A[It]}(u,v)$ can be used to compute
the Hilbert polynomial of the quotient ring $A/I^v$. In fact, we have
$$P_{A/I^v}(u) = P_A(u) - P_{A[It]}(u,v)$$
for $v$ large enough.

In particular, we can prove the following property of the function $e_{i}(M/I^kM)$ for any finitely generated graded $A$-module $M$ \cite{HPV}.

\begin{Theorem}[Herzog-Puthenpurakal-Verma, 2007]
The Hilbert coefficient $e_{i}(M/I^kM)$  as a function of $k$ is of
polynomial type of degree $\leq n-d+i$, where $d = \dim M/IM$. 
\end{Theorem}
 
Let  $V$ denote the blow-up of the subscheme of $\Proj A$ defined by 
$I$. It is known that $V$ can be embedded into a projective space by the 
linear system $(I^e)_c$ for any pair of positive integers $e,c$ with $c > 
de$ \cite{CH}.  Such embeddings often yield interesting rational 
varieties such as the Bordiga-White surfaces,
the Room surfaces and the Buchsbaum-Eisenbud varieties. Let $V_{c,e}$ denote the embedded variety.
The homogeneous coordinate ring of $V_{c,e}$ is the subalgebra $k[(I^e)_c]$ of $A$. 
It has been observed in \cite{STV} and \cite{CHTV} that $k[(I^e)_c]$ can be identified as the subalgebra of $A[It]$ along 
the diagonal $\{(cv,ev)|\ v \in \NN\}$ of $\NN^2$. 
Since $P_{A[It]}(cv,ev)$ is the Hilbert polynomial of $k[(I^e)_c]$, we 
may get uniform information 
on all such embeddings from $P_{A[It]}(u,v)$. 

\begin{Proposition} \label{embedded} Let $s = \dim A/0:I^\infty-1$. 
Assume that $c > de$. Then
$$\deg V_{c,e} = \displaystyle \sum_{i=0}^s \binom{s}{i}e_i(A[It])c^ie^{s-i}.$$ 
\end{Proposition}

If $I$ is generated by a $d$-sequence we have the following formula for the mixed multiplicities
$e_i(A[It])$, which displays a completely different behavior of 
$e_i(A[It])$ than that of $e_i(\mm|I)$ (see Theorem \ref{HTR}). 

\begin{Theorem}[Hoang-Trung, 2003] \label{d-sequence}
Let $I$ be an ideal generated by a homogeneous $d$-sequence 
$f_1,\ldots,f_r$ with $\deg f_j= d_j$ and $d_1\le\ldots\le d_r$. Let 
$I_q=(f_1,\ldots,f_{q-1}):f_q$ for $q=1,\ldots,r$. Set 
\begin{align*}
s & := \dim A/I_1-1,\\
m&  := \max\{q|\ \dim A/I_q + q-2 = s\}.
\end{align*}
Then $\deg P_{A[It]}(u,v) = s$ and
$$e_i(A[It]) =  
\sum_{q=1}^{\min\{m,s-i+1\}}(-1)^{s-q-i+1}e(A/I_q)\sum_{j_1+\ldots +j_q=s-q-i+1}d_1^{j_1}\ldots d_q^{j_q}$$
for $i = 0,\ldots,s$.
\end{Theorem}

Now we will apply Theorem \ref{d-sequence} to compute the mixed multiplicities of Rees algebras of complete intersections and of determinantal ideals.  

\begin{Corollary} \label{regular} 
Let $f_1,\ldots f_r$ be a homogeneous 
regular sequence with $\deg f_1=d_1 \le\ldots\le \deg f_r = d_r$ and $I 
= (f_1,\ldots,f_r)$. Set $s = \dim A-1$. Then $\deg P_{A[It]}(u,v)  = 
s$ and 
\begin{align*}
& e_i(A[It]) =\\
& \sum_{q=1}^{\min\{u,v-i+1\}} (-1)^{s-q-i+1} 
e(A)\sum_{j_1+\ldots +j_q=s-q-i+1}d_1^{j_1+1}\ldots d_{q-1}^{j_{q-1}+1}d_q^{j_q}
\end{align*}
for $i = 0,\ldots,s$.
\end{Corollary}

\begin{Corollary} \label{minor} 
Let $A = k[X]$, where $X$ is a 
$(r-1)\times r$ matrix of indeterminates. Let $I$ be the ideal of the maximal 
minors of $X$ in $A$. Set $s = (r-1)\times r-1$. Then $\deg 
P_{A[It]}(u,v)  = s$ and 
$$e_i(A[It]) = \sum_{q=1}^{\min\{u,v-i+1\}} (-1)^{s-q-i+1}\binom{r -1}{q-1} \binom{s-i}{q-1} r^{s-q-i+1}$$
for $i = 0,\ldots,s$.
\end{Corollary}

Hoang \cite{Ho2} also computed $e_i(A[It])$ in the case $I$ is the defining ideal of a rational normal curve.

\end{document}